\documentclass[a4paper,reqno,10pt]{amsart}

\usepackage{bbm} %bold numbers
\usepackage[T1]{fontenc} %allows for hyphenation of words with accents
\usepackage{newtxtext,newtxmath}
\usepackage{comment}
\usepackage{tikz-cd}

\theoremstyle{plain}
\newtheorem*{Theorem*}{Theorem} %no numbering for Theorem*

\theoremstyle{plain}
\newtheorem{Theorem}{Theorem}
\newtheorem{Lemma}[Theorem]{Lemma}
\newtheorem{Corollary}[Theorem]{Corollary}
\newtheorem{Proposition}[Theorem]{Proposition}

\theoremstyle{remark}
\newtheorem{Empty}[Theorem]{}

\newtheorem*{Claim*}{Claim}
\theoremstyle{exercise}

\numberwithin{Theorem}{section}
\numberwithin{Exercise}{section}

\newcommand{\R}{\mathbb{R}} %real numbers

\newcommand{\ind}{\mathbbm{1}} %indicatrix function

\newcommand{\bbE}{\mathbb{E}}

\newcommand{\bbP}{\mathbb{P}}

\newcommand{\calB}{\mathscr{B}}
\newcommand{\calC}{\mathscr{C}}
\newcommand{\calD}{\mathscr{D}}

\newcommand{\calF}{\mathscr{F}}
\newcommand{\calG}{\mathscr{G}}
\newcommand{\calH}{\mathscr{H}}

\newcommand{\calK}{\mathscr{K}}

\newcommand{\calP}{\mathscr{P}}

\newcommand{\calS}{\mathscr{S}}

\DeclareMathOperator{\rmdiam}{\operatorname{diam}} %diameter

\DeclareMathOperator{\rmspt}{\mathrm{spt}} %support
\newcommand{\rmtr}{\operatorname{tr}} %trace

\def\XXint#1#2#3{{%
\setbox0=\hbox{$#1{#2#3}{\int}$}
\vcenter{\hbox{$#2#3$}}\kern-.5\wd0}}

\newcommand{\la}{\langle}
\newcommand{\ra}{\rangle}

\renewcommand{\leq}{\leqslant}
\renewcommand{\geq}{\geqslant}
\renewcommand{\subset}{\subseteq}
\renewcommand{\supset}{\supseteq}

\newcommand{\niceBV}{\mathcal{BV}}
\newcommand{\fracdiv}{\mathfrak{div}}
\newcommand{\symdif}{\mathbin{\triangle}}

\begin{document}

%=================
% TITLE AND AUTHOR
%=================

%\titleAT %see above / page de garde d'un livre

\title{A regularity property of fractional Brownian sheets}

\author[Ph. Bouafia]{Philippe Bouafia}

\address{F\'ed\'eration de Math\'ematiques FR3487 \\
  CentraleSup\'elec \\
  3 rue Joliot Curie \\
  91190 Gif-sur-Yvette
}

\email{philippe.bouafia@centralesupelec.fr}

\author[Th. De Pauw]{Thierry De Pauw}

\address{Institute for Theoretical Sciences / School of Science, Westlake University\\
No. 600, Dunyu Road, Xihu District, Hangzhou, Zhejiang, 310030, China}

\email{thierry.depauw@westlake.edu.cn}

%\address{Universit\'e Paris Diderot\\ 
%  Sorbonne Universit\'e\\
%  CNRS\\ 
%  Institut de Math\'ematiques de Jussieu -- Paris Rive Gauche, IMJ-PRG\\
%  F-75013, Paris\\
%  France}
%\email{thierry.de-pauw@imj-prg.fr}



%\thanks{The first author was partially supported by the Science and Technology Commission of Shanghai (No. 18dz2271000).}

%=========
% ABSTRACT
%=========

%\begin{abstract}
%  We construct a predual of the space $BV_\infty(I^d)$ of bounded
%  functions of bounded variation over the square $I^d = [0,
%    1]^d$. Elements of this predual turn out to be in correspondence
%  with charges, an object that originates from nonabsolute integration
%  theories. We provide this space with a Faber-Schauder type basis,
%  which allows to derive chargeability and non-chargeability
%  criteria. As an application, we are able to show that the sample
%  paths of the Brownian sheet are almost surely not chargeable.
%\end{abstract}

\begin{abstract}
  A function $f$ defined on $[0, 1]^d$ is called strongly chargeable if
  there is a continuous vector-field $\varv$ such that $f(x_1, \dots,
  x_d)$ equals the flux of $\varv$ through the rectangle $[0, x_1]
  \times \cdots \times [0, x_d]$ for all $(x_1, \dots, x_d) \in
  [0, 1]^d$. In other words, $f$ is the primitive of the divergence of a
  continuous vector-field.
  We prove that the sample paths of the Brownian sheet with $d \geq 2$
  parameters are almost surely not strongly chargeable. On the
  other hand, those of the fractional Brownian sheet of Hurst parameter
  $(H_1, \dots, H_d)$ are shown to be almost surely strongly
  chargeable whenever
  \[
  \frac{H_1 + \cdots + H_d}{d} > \frac{d - 1}{d}.
  \]
\end{abstract}

\maketitle

\tableofcontents

\section{Introduction}

In order to motivate our results regarding multidimensional Brownian sheets we start with a few remarks about the 1-dimensional Brownian motion.
We recall that a Gaussian space is an infinite dimensional separable Hilbert space $E \subset L^2(\Omega,\calB,\bbP)$ containing only centered Gaussian variables where $(\Omega,\calB,\bbP)$ is a large enough probability space. There then exists a Gaussian noise, i.e.\ a Hilbertian isomorphism $G \colon L^2([0,1]) \to E$ and setting $W_t = G(\ind_{[0,t]})$, $0 \leq t \leq 1$, is a way of defining the standard Brownian motion $(W_t)_{0 \leq t \leq 1}$. Here, $\ind_{[0,t]}$ denotes the indicator function of the interval $[0,t]$. It is now trivial to check that the covariance $\bbE(W_s W_t) = \min(s,t)$. By an application of Kolmogorov's continuity theorem, one may assume that the function $t \mapsto W_t(\omega)$ is continuous for all $\omega \in \Omega$.
\par 
We say that $K \subset [0,1]$ is a {\it figure} if $K = \bigcup_{i=1}^p [s_i,t_i]$ for some finitely many pairwise nonoverlapping intervals $[s_1,t_1],\ldots,[s_p,t_p]$ and we let $\calF([0,1])$ denote the set of figures. The increment of the Brownian motion $W$ on $K$ is defined as
\[
\Delta_W K = \sum_{i=1}^p (W_{t_i} - W_{s_i}) = G(\ind_K).
\]
Even though $G(\ind_A)$ is defined for every measurable subset $A$ of $[0,1]$, it is not the case that, almost surely the increment $\calF([0,1]) \to \R \colon K \mapsto \Delta_W K$ extends to a (signed) Borel measure on $[0,1]$. Therefore, it makes sense to study the pathwise regularity of $K \mapsto \Delta_W K$ as a function of figures.
\par 
The increment is readily finitely additive, i.e.\ $\Delta_W(K_1 \cup K_2) = \Delta_W(K_1) + \Delta_W(K_2)$ whenever $K_1$ and $K_2$ are nonoverlapping. Furthermore, one easily checks that, owing to the continuity of $t \mapsto W_t$ the increment possesses the following continuity property: If $(F_n)_n$ is a sequence of figures whose (Lebesgue) measure tends to zero and whose number of components is uniformly bounded then $\Delta_W(F_n) \to 0$. A finitely additive function $\mu \colon \calF([0,1]) \to \R$ satisfying this continuity property is called a {\it charge} and it is easy to show that the space of charges $CH([0,1])$ is isomorphic with $C_0([0,1])$ by means of associating with $f \in C_0([0,1])$ the charge $\mu_f \colon \bigcup_{i=1}^p [s_i,t_i] \mapsto \sum_{i=1}^p (f(t_i)-f(s_i))$. Thus, it appears that the domain of a charge $\mu_f \in CH([0,1])$ can be extended to include all functions of bounded variation, referring to the Lebesgue-Stieltjes integral,
\[
\mu_f \colon BV([0,1]) \to \R \colon u \mapsto \int f du.
\] 
In this functional analytic context we think of charges as the members of the dual of $BV([0,1])$ with respect to some appropriate topology.
\par 
Finally, we recall a second point of view on Brownian motion, namely the L\'evy-Ciesielski construction. We denote by $h_{n,k}$ the Haar function supported in the interval $[k2^{-n},(k+1)2^{-n}]$ so that the sequence $(h_{n,k})_{n,k}$ is a Hilbertian basis of $L^2([0,1])$. The Faber-Schauder basis $(f_{n,k})_{n,k}$ of $C_0([0,1])$ is then obtained as a sequence of primitives of the former, $f_{n,k}(t) = \int_0^t h_{n,k}$. It was introduced by Faber in \cite{Fabe}. One can define the Brownian motion $W$ by its decomposition in the Faber-Schauder basis with an independent sequence of Gaussian centered coefficients $(A_{n,k})_{n,k}$:
\[
W_t = \sum_{n,k} A_{n,k} f_{n,k}(t).
\]
The advantage of this point of view is that one can study the regularity of a continuous function $\sum_{n,k} A_{n,k}f_{n,k}$ according to the asymptotic behavior of its sequence of coefficients $(A_{n,k})_{n,k}$, for instance whether it is H\"older continuous \cite{Cies}. See \cite{CiesKerkRoyn} for applications to probability. We are now ready to make sense of the corresponding observations for multidimensional stochastic processes.
\par 
Here, $d \geq 2$ is an integer. Given $H = (H_1,\ldots,H_d) \in (0,1)^d$ we say that a Gaussian centered random process $(W^H_{t_1,\ldots,t_d})_{0 \leq t_i \leq 1}$ is a fractional Brownian sheet of Hurst multiparameter $H$ if the covariance
\[
\bbE(W^H_{s_1,\ldots,s_d}W^H_{t_1,\ldots,t_d}) = \prod_{i=1}^d \frac{s_i^{2H_i} + t_i^{2H_i} - |t_i-s_i|^{2H_i}}{2}.
\]
When $H=(1/2,\ldots,1/2)$ we recover the standard Brownian sheet and we sometimes simply write $(W_{t_1,\ldots,t_d})_{0 \leq t_i \leq 1}$ with no reference to $H$. In that case,
\[
\bbE(W_{s_1,\ldots,s_d}W_{t_1,\ldots,t_d}) = \prod_{i=1}^d \min(s_i,t_i).
\]
As in case $d=1$, $(W_{t_1,\ldots,t_d})_{0 \leq t_i \leq 1}$ can be defined by means of a $d$-dimensional Gaussian noise. 
According to Kolmogorov's continuity theorem, one may assume that every sample of the function  $(t_1,\ldots,t_d) \mapsto W^H_{t_1,\ldots,t_d}$ is continuous.
\par 
We say that $K \subset [0,1]^d$ is a {\it rectangle} if $K = \prod_{i=1}^d [s_i,t_i]$ is a Cartesian product of compact intervals and we define {\it figures} to be the unions of finitely many rectangles (which we can assume, without loss of generality, are pairwise nonoverlapping, i.e.\ the Lebesgue measure of their intersection vanishes). The set of figures is denoted $\calF([0,1]^d)$. We also let $C_0([0,1]^d)$ be the space of continuous functions defined on $[0,1]^d$ that vanish at $x=(x_1,\ldots,x_d)$ if at least one $x_i=0$. With $f \in C_0([0,1]^d)$ is associated its increment $\Delta_fK$ on a rectangle $K$, whose definition we recall now only when $d=2$:
\begin{equation*}
%\label{eq.intro.1}
\Delta_f([s_1,t_1] \times [s_2,t_2]) = f(t_1,t_2) - f(s_1,t_2) - f(t_1,s_2) + f(s_1,s_2).
\end{equation*}
Since $\Delta_f$ is finitely additive on the set of rectangles, it extends uniquely to a finitely additive function of figures, still denoted $\Delta_f$. In general, $\Delta_f$ does not extend to a (signed) Borel measure on $[0,1]^d$. An interesting question consists of making sense of extra regularity properties of $\Delta_{X}$ when $X_{t_1,\ldots,t_d}$ is a multidimensional stochastic process with almost sure continuous realizations. This is what do in this paper in case $X = W^H$.
\par 
A $d$-dimensional {\it charge} is a finitely additive function $\mu \colon \calF([0,1]^d) \to \R$ satisfying the following continuity property: If $(F_n)_n$ is a sequence of figures such that $|F_n| \to 0$ and $\sup_n \|F_n\| \to 0$ then $\mu(F_n) \to 0$. Here, $|F_n|$ is the Lebesgue measure of the figure $F_n$ and $\|F_n\|=\calH^{d-1}(\partial F_n)$ is the $(d-1)$-dimensional Hausdorff measure of its boundary. Note that, when $d=1$ the term $\calH^0(\partial F)$ equals twice the number of components of the figure $F$, making clear the analogy with the higher dimensional case. The space of $d$-dimensional charges is denoted $CH([0,1]^d)$. In fact, $\mu \in CH([0,1]^d)$ extends (additively and continuously) to a larger collection of sets $A$ than figures, called sets of finite perimeter, i.e.\ whose indicator function $\ind_A \in BV([0,1]^d)$. Here, $BV([0,1]^d)$ is the space of functions of bounded variation in the sense of De Giorgi, i.e.\ those functions $u \in L^1([0,1]^d)$ whose distributional gradient is a vector-valued measure $Du$ of finite variation. See Section 3 for relevant information about functions of bounded variation and charges.
\par 
Recall that $d \geq 2$. We prove that:
\begin{itemize}
\item The increment $\Delta_W$ of the sample paths of the Brownian sheet $W$ are almost surely not a charge;
\item The increment $\Delta_{W^H}$ of the sample paths of the fractional Brownian sheet are almost surely a charge provided $H_1 + \cdots + H_d > d-1$.
\end{itemize}
In fact, we prove more than this. In order to state our results, we need to introduce the notion of a {\it strong charge}. We consider the continuous linear embedding $T \colon L^d([0,1]^d) \to BV([0,1]^d)^*$ defined by $T_f(u) = \int fu$. We let $SCH([0,1]^d)$ be the closure of the range of $T$ and we call its members the $d$-dimensional strong charges. One can show that $SCH([0,1]^d)$ is a predual of $BV([0,1]^d)$, theorem \ref{thm:predual}. The strong charges are exactly the linear functionals $\alpha \colon BV([0,1]^d) \to \R$ associated with a continuous vector-field $\varv \colon [0,1]^d \to \R^d$ in the following way: $\alpha(u) = - \int \la \varv , Du \ra$ for each $u \in BV([0,1]^d)$, theorem \ref{thm:fracdivonto}. With each strong charge $\alpha$ one can associate a (unique) charge $\mu_\alpha$ by means of the formula $\mu_\alpha(A) = \alpha(\ind_A)$ but the converse is not true.
\par 
We identify a Schauder basis of $SCH([0,1]^d)$ and we establish useful criteria for corresponding random series to be strong charges. These apply to the Brownian and fractional Brownian sheet. We start by describing multidimensional Haar functions (see \ref{emp:multiHaar}). These are tensor products of their 1-dimensional analogues. Specifically, we let
\[
 g_{n, k, r} = 2^{nd/2} \sum_{\ell=0}^{2^d - 1} (A_d)_{r, \ell}
    \ind_{K_{n+1, 2^d k + \ell}}
\]
where the orthogonal $2^d$-dimensional square matrix $A_d$ is the Kronecker product of $\begin{pmatrix}1 & 1 \\ 1 & -1\end{pmatrix}$, $d$ times with itself, and the dyadic cubes $K_{n,j}$ have been numbered in an appropriate way where $n$ denotes their generation, i.e.\ $K_{n,j}$ has side length $2^{-n}$. One shows that the sequence $(g_{n,k,r})_{n,k,r}$ is a Hilbertian basis of $L^2([0,1]^d)$. Similarly to the 1-dimensional case, one then wants to consider ``primitives'' of these Haar functions in some sense, as did Faber. We introduce the strong charges $T_{g_{n,k,r}}$ (recall the embedding $T$ from the previous paragraph) as primitives of the Haar functions $g_{n,k,r}$ and we prove these constitute a Schauder basis of $SCH([0,1]^d)$, theorem \ref{thm:mainSchauder}. In fact, we show that if $\alpha \in SCH([0,1]^d)$ then
\[
\alpha = \sum_{n,k,r} \alpha(g_{n,k,r}) T_{g_{n,k,r}}.
\]
The strong charges $T_{g_{n,k,r}}$ behave somewhat like a wavelet basis in the sense that they have nonoverlapping supports, $\rmspt(T_{g_{n,k,r}}) = K_{n,k}$. Taking advantage of this, we are able to determine whether some series $\sum_{n,k,r} a_{n,k,r}T_{g_{n,k,r}}$ converge or not to a strong charge, according to the asymptotic behavior of their coefficients $(a_{n,k,r})_{n,k,r}$,  corollary \ref{cor:chargeability}.
\par 
Now let $f \in C_0([0,1]^d)$ and $\Delta_f$ be its increment. If $\alpha = \Delta_f$ were a strong charge then the coefficients in the above convergent series would be
\[
\Delta_f(g_{n,k,r}) = 2^{nd/2} \sum_{\ell=0}^{2^d - 1} (A_d)_{r, \ell}
    \Delta_f(\ind_{K_{n+1, 2^d k + \ell}}).
\]
In case $f$ is $W^H$, the increments $\Delta_{W^N}(K_{n,j})$ on dyadic cubes $K_{n,j}$ are random variables whose asymptotics as $(n,j)$ grows can be controlled. Together with the quantitative criteria evoked at the end of last paragraph, we are then able to establish that:
\begin{itemize}
\item The increment $\Delta_{W^H}$ of the sample paths of the fractional Brownian sheet $W^H$ are almost surely not a charge if $\bar{H} \leq \frac{d-1}{d}$, theorem \ref{thm:counterEx};
\item The increment $\Delta_{W^H}$ of the sample paths of the fractional Brownian sheet $W^H$ are almost surely a strong charge if $\bar{H} > \frac{d-1}{d}$, theorem \ref{thm:mainFrac};
\end{itemize}
where
\[
\bar{H} = \frac{H_1 + \ldots + H_d}{d}.
\]
\par 
Charges are natural integrators in non-absolute integration theories,
see for example \cite{Pfef4}. We already mentionned that strong
charges act on $BV$ functions. Thus, in case $\bar{H} >
\frac{d-1}{d}$, it makes sense to integrate $BV$ functions with
respect to the charge $\Delta_{W^H}$ representing the variations of
the fractional Brownian sheet. Such an integral is understood in a
pathwise sense. The possibility to extend this integral to a
full-fledged Young integral, where the integrand is allowed to be a
Hölder continuous function, will be investigated in \cite{Boua}.
\section{Notations}

Throughout this paper, $\R$ will denote the set of real
numbers. We will work in an ambient space whose dimension is an
integer $d \geq 1$, typically $(0, 1)^d$, $[0, 1]^d$, or $\R^d$. The
Euclidian norm of $x \in \R^d$ is denoted $|x|$.

The closure, the interior, and the topological boundary of a set $E
\subset \R^d$ will be denoted $\operatorname{cl} E$,
$\operatorname{int} E$, and $\partial E$, respectively. The indicator
function of $E$ is $\ind_E$. The symmetric difference of two sets
$E_1, E_2 \subset \R^d$ is written $E_1 \symdif E_2$.

Unless otherwise specified, the expressions ``measurable'', ``almost
all'', and ``almost everywhere'' tacitly refer to the Lebesgue
measure. The Lebesgue (outer) measure of a set $E \subset \R^d$ is
simply written $|E|$. Two subsets $E_1, E_2 \subset \R^d$ are said to
be almost disjoint whenever $|E_1 \cap E_2| = 0$. If $U \subset \R^d$
is a measurable set and $1 \leq p \leq \infty$, the corresponding Lebesgue spaces
are denoted $L^p(U)$. Here, $U$ is endowed with its Lebesgue
$\sigma$-algebra and the Lebesgue measure. The $L^p$ norm is written
$\|\cdot\|_p$. The notation $\|\cdot\|_\infty$ might also refer to the
supremum norm in the space of continuous functions. The integral of a
function $f$ with respect to the Lebesgue measure is simply written
$\int f$, with no mention of the Lebesgue measure. In case another measure
is used, it will be clear from the notation. The $(d-1)$-dimensional
Hausdorff measure (defined on Borel subsets of $\R^d$) is denoted
$\calH^{d-1}$ and the corresponding $L^p$ spaces are written $L^p(U;
\calH^{d-1})$, where $U$ is a Borel subset of $\R^d$.

The topological dual of a Banach space $X$ is $X^*$. Unless otherwise
specified, the operator norm of a continuous linear map $T$ between
normed spaces will be written $\|T\|$.

\section{Preliminaries on $BV$ functions, $BV$ sets, and charges}

%In this section, we will properly define the space $SCH(I^d)$ of
%divergences of continuous vector fields on $I^d$. For historical
%reasons, elements of $SCH(I^d)$ are called strong charges. The results
%of this paper rely heavily on functional analytic methods. In fact,
%the reason we restrict ourselves to studying strong charges on the
%bounded domain $I^d$ is that we want $SCH(I^d)$ to be a Banach
%space. The duality between strong charges and functions of bounded
%variation will also play an important role, so much so that in our
%presentation, $SCH(I^d)$ will be constructed as a predual of $BV(I^d)$,
%already embedded in $BV(I^d)^*$.

\begin{Empty}[$BV$ functions]
  We start by introducing all the necessary definitions and results
  concerning functions of bounded variation. For more insight about these,
  we refer to the book of Evans and Gariepy~\cite{EvanGari}.

  Let $U \subset \R^d$ be an open set. The variation of a Lebesgue
  integrable function $u \colon U \to \R$ over an open subset $V
  \subset U$ is the quantity
  \begin{equation}
    \label{def:var}
    \|Du\|(V) = \sup \left\{ \int_V u \operatorname{div} \varv : \varv
    \in C^1_c(V; \R^d) \text{ and } |\varv(x)| \leq 1 \text{ for all }
    x \in V \right\}
  \end{equation}
  where $C^1_c(V; \R^d)$ denotes the space of continuously
  differentiable compactly supported vector fields on $V$.

  The function $u$ is said to be of bounded variation whenever
  $\|Du\|(U) < \infty$. In this case, the vector-valued
  Riesz representation theorem can be used to prove that the
  distributional gradient of $u$ is an $\R^d$-valued Borel measure
  denoted $Du$. Its total variation measure is denoted $\|Du\|$; it is
  a finite Borel measure whose values on open subsets $V$ of $U$ is
  given by formula~\eqref{def:var}.

  The set of (equivalence classes of) functions of bounded variation
  on $U$ is denoted $BV(U)$. It is a Banach space under the norm
  $\|u\|_{BV} = \|u\|_1 + \|Du\|(U)$. In the sequel, the following
  results will be used.
  \begin{Theorem*}[Compactness theorem, {\cite[5.2.3, Theorem 4]{EvanGari} and \cite[5.2.1, Theorem 1]{EvanGari}}]
    Let $U \subset \R^d$ be a bounded Lipschitz open set and $(u_n)$
    be a bounded sequence in $BV(U)$. There is a subsequence
    $(u_{n_k})$ and a function $u \in BV(U)$ such that $u_{n_k} \to u$
    in $L^1(U)$. Furthermore, $\|Du\|(U) \leq \liminf
    \|Du_{n_k}\|(U)$.
  \end{Theorem*}
  One of the consequences of the compactness theorem is that the
  closed unit ball of $BV(U)$ is compact when given the
  $L^1$-topology. This is a strong indication that $BV(U)$ is a dual
  Banach space. Indeed, this will be proven in
  Theorem~\ref{thm:predual} (for the case $U = (0, 1)^d$, but the proof
  applies to any bounded Lipschitz open set as well), see also \cite[Remark 3.12]{AMBROSIO.FUSCO.PALLARA} for another point of view.
  \begin{Theorem*}[Sobolev-Poincaré inequality, {\cite[Theorem 5.11.1]{Ziem}}]
    Let $U \subset \R^d$ be a connected bounded Lipschitz open set and $\gamma
    \in BV(U)^*$ be a continuous linear functional such that
    $\gamma(\ind_U) = 1$. There is a constant $C = C(U, \gamma) \geq
    0$ such that for all $u \in BV(U)$,
    \[
    \|u - \gamma(u)\|_{d / (d-1)} \leq C \|Du\|(U)
    \]
  \end{Theorem*}
  In the above statement as well as in the remaining part of this paper, when $d=1$ we agree that $d/(d-1) = \infty$.
  We will apply the Sobolev-Poincar\'e inequality to the case where
  \begin{equation}
    \label{eq:average}
  \gamma(u) = \frac{1}{|U|} \int_U u \text{ for all } u \in
  BV(U)
  \end{equation}
  and more specifically to domains that are $d$-dimensional open
  cubes, that is, sets of the form $U = \prod_{i=1}^d (a_i, b_i)$,
  where $|b_1 - a_1| = \cdots = |b_d - a_d| > 0$. A simple scaling
  argument shows that, in case of the averaging functional $\gamma$ defined
  in~\eqref{eq:average}, the Poincaré constant is the same for all
  $d$-dimensional open cubes. We will denote this constant by $C_P$.

  \begin{Theorem*}[Trace theorem, {\cite[5.3, Theorem 1]{EvanGari}}]
    If $U \subset \R^d$ is a bounded Lipschitz open set, there is a
    continuous linear operator $\rmtr \colon BV(U) \to L^1(\partial U;
    \calH^{d-1})$ such that for all $u \in BV(U)$ and $\varv \in C^1(\R^d;
    \R^d)$,
    \[
    \int_U u \operatorname{div} \varv = \int_{\partial U} \rmtr(u)
    \varv \cdot n_{U} d \calH^{d-1} - \int \varv \cdot dDu
    \]
    where $n_{U}$ denotes the normal outer unit vector field defined
    $\calH^{d-1}$-almost everywhere on $\partial U$.
  \end{Theorem*}

  A useful corollary of the trace theorem is the following result.
  
  \begin{Theorem*}[Extension theorem, {\cite[5.4, Theorem 1]{EvanGari}}]
    If $U \subset \R^d$ is a bounded Lipschitz open set and $u \in
    BV(U)$, define the function
    \[
    Eu \colon x \mapsto \begin{cases} u(x) & \text{if } x \in U \\
      0 & \text{if } x \not\in U
    \end{cases}
    \]
    Then $Eu \in BV(\R^d)$ and
    \[
    \|D(Eu)\|(\R^d) = \|Du\|(U) + \int_{\partial U} |\rmtr u| d \calH^{d-1}
    \]
  \end{Theorem*}
\end{Empty}

\begin{Empty}[$BV$ sets and charges]
  \label{e:charges}
  The \emph{perimeter} of a measurable subset $A$ of $\R^d$ is the
  extended real number $\|A\|= \|D\ind_A\|(\R^d)$. Usually, $A$ is
  said to be a set of finite perimeter (or a Caccioppoli set)
  whenever $\|A\| < \infty$. However, in this paper, we will rather
  follow Pfeffer's terminology \cite{Pfef2}: we will say that $A$ is a
  \emph{$BV$-set} whenever $A$ is bounded, measurable and $\|A\| <
  \infty$.

  We let $\niceBV(A)$ be the set of $BV$-subsets of $A$.
  %It can be proven that $\niceBV(A)$ is an algebra of sets, as, for
  %all $B, B_1, B_2 \in \niceBV(A)$, one has
  %\begin{equation}
  %  \label{eq:perimeters}
  %  \|A \setminus B\| \leq \|A\| + \|B\| \text{ and } \max(\| B_1 \cup
  %  B_2 \|, \|B_1 \cap B_2\|) \leq \|B_1\| + \|B_2\|
  %\end{equation}%
  It is endowed with the following notion of
  convergence: a sequence $(B_n)$ is said to \emph{$BV$-converge} to
  $B$ whenever
  \[
  \sup \|B_n\| < \infty \text{ and } \lim |B_n \symdif B| = 0.
  \]
  There exists a topology on $\niceBV(A)$ that is compatible with this
  notion of convergent sequences (see \cite[p. 33 and p. 42]{Pfef2}), but we shall not use it in the present paper.
  %Using~\eqref{eq:perimeters}, it may be
  %proven that the set-theoretic operations are continuous with respect
  %to $BV$-convergence. For example, the following will be used in the
  %paper: if two sequences $(B_n)$ and $(C_n)$ respectively
  %$BV$-converge to $B$ and $C$, then $(B_n \setminus C_n)$
  %$BV$-converge to $B \setminus C$.

  A \emph{charge on $A$} is a function $\mu \colon \niceBV(A) \to
  \R$ that satisfies the following proporties:
  \begin{enumerate}
  \item[(A)] Finite additivity: $\mu(B_1 \cup B_2) = \mu(B_1) +
    \mu(B_2)$ whenever $B_1, B_2 \in \niceBV(A)$ are almost
    disjoint;
  \item[(B)] Continuity with respect to $BV$-convergence: if a
    sequence $(B_n)$ $BV$-converges to $B$, then $\mu(B_n) \to
    \mu(B)$.
  \end{enumerate}
  We observe that a charge necessarily vanishes on negligible
  sets. This can be seen as a consequence of either (A) or (B).  The
  linear space of charges on $A$ is denoted $CH(A)$. It is worth mentioning that
  our notation $\mu$ is not meant to suggest that charges are measures -- indeed, 
  some are not. However, absolutely continuous measures are charges. 

  The structure of $1$-dimensional $BV$-sets is strikingly simple and
  this allows for an easy description of charges in case $A = [0,
    1]$. Indeed, the elements of $\niceBV([0, 1])$ are, up to
  negligible sets, the disjoint unions of compact intervals. From (B), the
  function $\varv \colon [0, 1] \to \R$ defined by $\varv(x) = \mu([0,
    x])$ is continuous and vanishes at $0$. Reciprocally, to any
  function $\varv$ belonging to the space $C_0([0, 1])$ of continuous
  functions on $[0, 1]$ vanishing at $0$, we associate the charge
  $\Delta_\varv$ that maps any disjoint union of compact intervals to
  \[
  \Delta_{\varv} \colon \bigcup_{i=1}^p [a_i, b_i] \mapsto
  \sum_{i=1}^p \varv(b_i) - \varv(a_i)
  \]
  and such that $\Delta_v(B_1) = \Delta_v(B_2)$ whenever $B_1, B_2 \in
  \niceBV([0, 1])$ are almost disjoint, {\it i.e.} $|B_1 \symdif B_2| = 0$. Thus, $\varv
  \mapsto \Delta_{\varv}$ is a bijection from
  $C_0([0, 1])$ to $CH([0, 1])$.
\end{Empty}

\begin{Empty}[The Banach space {$C_0([0, 1]^d)$}]
  The multidimensional generalization of the space $C_0([0, 1])$ is
  the space $C_0([0, 1]^d)$ of continuous functions on $[0, 1]^d$ that
  vanish on the coordinate hyperfacets
  \begin{equation}
    \label{eq:axes}
    \bigcup_{i=1}^d \left\{ (x_1, \dots, x_d) \in [0, 1]^d : x_i = 0
    \right\}.
  \end{equation}
  We equip this space with the maximum norm $\| \cdot \|_\infty$.
\end{Empty}

\begin{Empty}[Chargeability]
  Among the subsets of $[0, 1]^d$, we shall consider some that are more
  regular than $BV$ sets. We describe these here.

  A \emph{dyadic cube} is a set of the type $\prod_{i=1}^d
  \left[2^{-n} k_i, 2^{-n} (k_i+1) \right]$, where $n \geq 0$ and $0
  \leq k_1, \dots, k_d \leq 2^n-1$ are integers. Such a dyadic cube
  has side length $2^{-n}$ and we will say that it is \emph{of
    generation} $n$. Thus, in our terminology, dyadic cubes are subsets of $[0,1]^d$.

  A \emph{rectangle} is a set with non-empty interior of the form
  $\prod_{i=1}^d [a_i, b_i]$. We will also be interested in
  \emph{(rectangular) figures} in $[0, 1]^d$, \textit{i.e} subsets of
  $[0, 1]^d$ that can be written as finite unions of rectangles. The
  collection of such sets is denoted $\calF([0, 1]^d)$. Note that each
  $F \in \calF([0, 1]^d)$ is a $BV$-set and that $\|F\| = \calH^{d-1}(\partial F)$.

  We can now define the increments of a function $f \in C_0([0, 1]^d)$
  on a rectangle $\prod_{i=1}^d [a_i, b_i] \subset [0, 1]^d$ by means
  of the formula
  \[
  \Delta_f\left( \prod_{i=1}^d [a_i, b_i] \right) = \sum_{(c_i) \in
    \prod_{i=1}^d \{a_i, b_i\}} (-1)^{\delta_{a_1, c_1}} \ldots
  (-1)^{\delta_{a_d, c_d}} f(c_1, \dots, c_d)
  \]
  where for any reals $t, t'$, the number $\delta_{t, t'}$ is $1$ or
  $0$ according to whether $t = t'$ or not. In the $2$-dimensional
  case, we recover the well-known rectangular increment
  $\Delta_f([a_1, b_1] \times [a_2, b_2]) = f(b_1, b_2) - f(a_1, b_2)
  - f(b_1,a_2) + f(a_1, a_2)$.

  One checks that if a rectangle $K$ can be split as the
  union of two almost disjoint rectangles $K = K_1 \cup K_2$, than
  $\Delta_f(K) = \Delta_f(K_1) + \Delta_f(K_2)$. Based on this
  observation and the fact that any figure can be written as a finite
  union of pairwise almost disjoint rectangles, $\Delta_f$ has a
  unique extension to $\calF([0, 1]^d)$, that satisfies the finite
  additivity property (A) of the paragraph~\ref{e:charges}, restricted
  to the subcollection $\calF([0, 1]^d)$ of $\niceBV([0, 1]^d)$. This
  extension is still denoted $\Delta_f$.

  Now, we say that the function $f$ is \emph{chargeable} whenever
  $\Delta_f \colon \calF([0, 1]^d) \to \R$ has an extension to
  $\niceBV([0, 1]^d)$ that is a charge. We state below an
  approximation theorem of De Giorgi (see \cite[Proposition 1.10.3]{Pfef2} for a proof), that implies that this extension is
  necessarily unique. In intuitive terms, the chargeability of $f$
  allows to make sense of increments of $f$ over arbitrary
  $BV$-sets. The discussion in~\ref{e:charges} shows that all
  continuous functions on $[0, 1]$ vanishing at $0$ are chargeable. In
  fact, chargeability is a regularity property that differs from continuity only in
  dimension $\geq 2$. We will prove later that the sample paths of the
  Brownian sheet are almost surely not chargeable.

  Functions that are chargeable can be thought of as being primitives
  of charges, by integration on rectangles. This statement is made
  clear by the elementary Proposition~\ref{prop:primcharge}.  The
  class of $BV$-subsets of $[0, 1]^d$ has better properties than that
  of rectangular figures. It is invariant under biLipschitz
  transformations and its definition does not rely on a specific
  choice of a basis in $\R^d$. We uphold the thesis that, whenever a
  function $f$ is chargeable, the charge $\Delta_f$ is a more
  fundamental object than $f$ itself.
\end{Empty}

\begin{Theorem}[De Giorgi approximation]
  There is a constant $C \geq 0$, depending only on $d$, such
  that, for any $BV$-subset $B \subset [0, 1]^d$, there exists a
  sequence $(F_n)$ of figures in $[0,1]^d$ such that
  \[
  \sup \|F_n\| \leq C \|B\| \text{ and } \lim |F_n \symdif B| = 0.
  \]
  In particular, $(F_n)$ $BV$-converges to $B$.
\end{Theorem}

\begin{Proposition}
  \label{prop:primcharge}
  A function $f \in C_0([0, 1]^d)$ is chargeable if and only if there
  exists a charge $\mu$ on $[0, 1]^d$ such that
  \begin{equation}
    \label{eq:primcharge}
  f(x_1, \dots, x_d) = \mu\left( \prod_{i=1}^d [0, x_i] \right), \qquad 0 \leq x_1, \dots, x_d \leq 1 .
  \end{equation}
  In this case, one has $\mu = \Delta_f$.
\end{Proposition}

\begin{proof}
  If $f$ is chargeable, then $\Delta_f$ is a charge and $f(x_1, \dots,
  x_d) = \Delta_f([0, x_1] \times \cdots \times [0, x_d])$, as $f$
  vanishes on~\eqref{eq:axes}. Conversely, suppose the existence of a
  charge $\mu$ such that~\eqref{eq:primcharge} holds and consider a
  rectangle $K = \prod_{i=1}^d [a_i, b_i]$. Then, by the finite
  additivity of $\mu$, one obtains
  \[
  f(x_1, \dots, x_{d-1}, b_d) - f(x_1, \dots, x_{d-1}, a_d) =
  \mu\left( \prod_{i=1}^{d-1} [0, x_i] \times [a_d, b_d]\right)
  \]
  for any $0 \leq x_1, \dots, x_{d-1} \leq 1$. Repeating this process
  one step further, one gets
  \begin{multline*}
  \sum_{c_{d-1}, c_d} (-1)^{\delta_{a_{d-1}, c_{d-1}}}
  (-1)^{\delta_{a_d, c_d}} f(x_1, \dots, x_{d-2}, c_{d-1}, c_d) = \\
  \mu \left( \prod_{i=1}^{d-2} [0, x_i] \times [a_{d-1},
    b_{d-1}] \times [a_d, b_d]\right)
  \end{multline*}
  where $c_{d-1}$ ranges over $\{a_{d-1}, b_{d-1}\}$ and $c_d$ over
  $\{a_d, b_d\}$.  Continuing further, we obtain $\Delta_f(K) =
  \mu(K)$. By finite additivity of both $\Delta_f$ and $\mu$, we
  deduce that $\mu = \Delta_f$ on $\calF([0, 1]^d)$, which yields the
  result.
\end{proof}

\begin{Empty}[The space {$BV([0, 1]^d)$}]
  \label{e:defBVclosed}
  We let $BV([0, 1]^d)$ be the subspace of $BV(\R^d)$ that consists of
  those functions $u$ such that $u = 0$ almost everywhere on $\R^d
  \setminus [0, 1]^d$, equipped with the norm inherited from $BV(\R^d)$.

  In fact, the extension operator $E \colon BV((0, 1)^d) \to BV([0,
    1]^d)$ provides an isomorphim, whose reciprocal is the restriction
  operator. A measurable function $u \colon \R^d \to \R$ belongs to
  $BV([0, 1]^d)$ if and only if $u = 0$ almost everywhere outside $[0,
    1]^d$ and the restriction of $u$ to $(0, 1)^d$ belongs to
  $BV((0,1)^d)$. This means that the norms $\| \cdot \|_{BV(\R^d)}$ and
  \[
  u \mapsto \int_{(0, 1)^d} |u| + \|Du\|((0, 1)^d)
  \]
  are equivalent on $BV([0,1]^d)$.

  We claim that $u \mapsto \|Du\|(\R^d)$ is yet another norm on
  $BV([0, 1]^d)$ that is equivalent to the two norms above. To prove
  this claim, let $\gamma \colon BV((0, 1)^d) \to \R$ be the map
  defined by
  \[
  \gamma(\tilde{u}) = \frac{1}{2d}
  \int_{\partial(0, 1)^d} \rmtr \tilde{u} \, d\calH^{d-1}, \qquad
  \tilde{u} \in BV((0, 1)^d) .
  \]
  It is continuous, by the continuity of $\rmtr$. By the
  Sobolev-Poincaré inequality, there is a constant $C$ such that
  $\|\tilde{u} - \gamma(\tilde{u})\|_{d/(d-1)} \leq C
  \|D\tilde{u}\|((0, 1)^d)$ for all $\tilde{u}$. Furthermore, one has,
  by the extension theorem,
  \begin{equation}
    \label{eq:varUeq}
    \|Du\|(\R^d) = \|D\tilde{u}\|((0, 1)^d) + \int_{\partial (0, 1)^d} | \operatorname{tr} \tilde{u} | d \calH^{d-1}
  \end{equation}
  where $\tilde{u} \in BV((0, 1)^d)$ denotes here the restriction of
  $u$.  One notices that
  \begin{align*}
  \|u\|_1 = \|\tilde{u}\|_1 & \leq \|\tilde{u} - \gamma(\tilde{u})\|_1
  + \frac{1}{2d} \|\rmtr \tilde{u}\|_1 \\ & \leq \|\tilde{u} -
  \gamma(\tilde{u}) \|_{d/(d-1)} + \frac{1}{2d} \|\rmtr \tilde{u}\|_1
  \\ & \leq \max \left(C, \frac{1}{2d} \right) \left(
  \|D\tilde{u}\|((0, 1)^d) + \|\rmtr \tilde{u}\|_1 \right)
  \end{align*}
  Thus by~\eqref{eq:varUeq}, one obtains $\|u\|_{BV} = \|u\|_1 +
  \|Du\|(\R^d) \leq C' \|Du\|(\R^d)$ for some constant $C' > 0$. The
  upper bound $\|Du\|(\R^d) \leq \|u\|_{BV}$ is trivial.

  For the definition of strong charge functionals in the next section,
  it will be slightly more convenient to work with $BV([0, 1]^d)$
  rather than the isometric space $BV((0, 1)^d)$. One could similarly
  define a Banach space $BV(A)$, where $A$ is a $BV$-set, and develop
  a theory of strong charge functionals on $A$.
\end{Empty}

\section{Strong charge functionals and strong charges}

\begin{Empty}[The space of strong charge functionals]
  To each function $f \in L^d([0, 1]^d)$, we associate the functional
  \[
  T_f \colon u \mapsto \int_{[0, 1]^d} fu
  \]
  defined on $BV([0, 1]^d)$. It is continuous. Indeed, denoting $\bar{u}$
  the average value of $u$ on $[0, 1]^d$, we have, by the Hölder and
  Sobolev-Poincaré inequalities,
  \begin{align*}
  |T_f(u)| & = \left|\int_{[0, 1]^d} f(u - \bar{u}) + \left(
  \int_{[0, 1]^d} f \right) \bar{u} \right| \\
  & \leq \|f\|_d \left(\int_{(0, 1)^d} |u - \bar{u}|^{d/(d-1)}\right)^{(d-1)/d} +  \|f\|_d \|u\|_1 \\
  & \leq \|f\|_d \left( C_P \|Du\|( (0,1)^d ) + \|u\|_1 \right) \\
  & \leq \max(C_P, 1) \|f\|_d \|u\|_{BV}.
  \end{align*}
  Furthermore, this computation shows that the map $T \colon L^d([0, 1]^d)
  \to BV([0, 1]^d)^*$ sending $f$ to $T_f$ is continuous.

  We define the Banach space $SCH([0, 1]^d)$ as the closure of the range of
  the operator $T$ in $BV([0, 1]^d)^*$. Elements thereof are called \emph{strong
    charge functionals}. The choice of the terminology and the link with 
    charges defined in the previous section will be explained in \ref{strong_charges}.
    With this definition, $SCH([0, 1]^d)$ turns out to be a
  predual of $BV([0, 1]^d)$, as we prove now.
\end{Empty}

\begin{Theorem}[Duality theorem]
  \label{thm:predual}
  The canonical map
  \[
  \Upsilon \colon BV([0, 1]^d) \to SCH([0, 1]^d)^*
  \]
  (that sends $u$ to the functional $\alpha \mapsto \alpha(u)$) is an
  isomorphism of Banach spaces.
\end{Theorem}

\begin{proof}
  The canonical map $\Upsilon$ is the composition of the evaluation
  map $BV([0, 1]^d) \to BV([0, 1]^d)^{**}$ with the adjoint map of the injection
  map $SCH([0, 1]^d) \to BV([0, 1]^d)^*$, and, therefore, is continuous. Next, we
  prove separately that $\Upsilon$ is one-to-one and onto.

  If $u$ belongs to the kernel of $\Upsilon$, then
  \[
  T_f(u) = \int_{[0, 1]^d} fu = 0
  \]
  for all $f \in L^d([0, 1]^d)$. From this it can be deduced that $u
  = 0$ a.e on $[0, 1]^d$, hence $u = 0$ in $BV([0, 1]^d)$. Therefore,
  $\Upsilon$ is an injection.

  Now let us take $\gamma \in SCH([0, 1]^d)^*$. As $\gamma \circ T$
  belongs to the dual of $L^d([0, 1]^d)$, there is a function $u \in
  L^{d/(d-1)}([0, 1]^d)$ such that
  \[
  \gamma(T_f) = \int_{[0, 1]^d} fu
  \]
  for all $f \in L^d([0, 1]^d)$. We extend $u$ to $\R^d$ by zero, and
  we wish to prove that $u$ belongs to $BV([0, 1]^d)$. To this end, we consider a vector field $\varv \in C^1_c(\R^d; \R^d)$
  such that $|\varv| \leq 1$ on $\R^d$. Call $g$ the restriction of
  $\operatorname{div} \varv$ to $[0, 1]^d$. First, note that
  \[
  T_{g}(\varphi) = \int_{[0, 1]^d} \varphi
  \operatorname{div} \varv = \int_{\R^d} \varphi \operatorname{div} \varv \leq \|D\varphi\|(\R^d) \leq
  \|\varphi\|_{BV}
  \]
  for all $\varphi \in BV([0, 1]^d)$. This establishes that
  $\|T_{g}\| \leq 1$. Then we
  observe that
  \[
  \int_{[0, 1]^d} u \operatorname{div} \varv = \gamma(T_{g}) \leq
  \|\gamma\|.
  \]
  As $\varv$ is arbitrary, this proves that $\|Du\|(\R^d) \leq
  \|\gamma\| < \infty$ and so $u \in BV([0, 1]^d)$.

  The continuous maps $\gamma$ and $\Upsilon(u)$ coincide on $T(L^d([0,
    1]^d))$, a dense subspace of $SCH([0, 1]^d)$. On this
  account, we infer that $\gamma = \Upsilon(u)$. So, $\Upsilon$ is onto.
  Finally, we apply the open mapping theorem to conclude that
  $\Upsilon^{-1}$ is continuous as well.
\end{proof}

The next proposition characterizes weak* convergence of sequences in
$BV([0, 1]^d)$ (weak* convergence with respect to the duality between $SCH([0, 1]^d)$ and $BV([0,
  1]^d)$). This is the same convergence that appears in the $BV$
compactness theorem.

\begin{Proposition}
  \label{prop:weak*conv}
  A sequence $(u_n)$ in $BV([0, 1]^d)$ weak* converges to $u$ if and
  only if it is bounded and $u_n \to u$ in $L^1$.
\end{Proposition}

\begin{proof}
  Of course, a weak* convergent sequence $(u_n)$ is bounded by the
  uniform boundedness principle. To prove that $u_n \to u$ in $L^1([0,
    1]^d)$, it suffices to remark that any subsequence of $(u_n)$ has
  a subsequence converging to $u$ in $L^1(\R^d)$. Indeed, by the
  compactness theorem (that we may apply to the functions $u_n$
  restricted to a bounded Lipschitz open neighborhood of $[0, 1]^d$),
  it is possible to extract, from any subsequence $(u_{n_k})$ of
  $(u_n)$, a subsequence still denoted $(u_{n_k})$ that
  $L^1$-converges to some $\varv \in BV([0, 1]^d)$. For any $f \in
  L^\infty([0, 1]^d)$, we have $T_f(u_{n_k}) \to T_f(\varv)$, whereas
  we also have $T_f(u_{n_k}) \to T_f(u)$ by weak* convergence. Hence
  $T_f(u) = T_f(\varv)$. As $f \in L^\infty([0, 1]^d)$ is arbitrary,
  we deduce that $u = \varv$ and this finishes the first part of the
  proof.

  Conversely, we observe that the space $L^\infty([0, 1]^d)$ is dense
  in $L^d([0, 1]^d)$, and therefore $T(L^\infty([0, 1]^d))$ is dense
  in $SCH([0, 1]^d)$. Owing to the boundedness of $(u_n)$, it is
  sufficient to remark that $T_f(u_n) \to T_f(u)$ for any $f \in
  L^\infty([0, 1]^d)$.
\end{proof}

\begin{Empty}[The operator $\fracdiv$]
  To each continuous vector field $\varv \in C([0, 1]^d; \R^d)$, we
  associate the functional $\fracdiv\, \varv : BV([0, 1]^d) \to \R$ defined by
  \[
  (\fracdiv \, \varv)(u) =
  - \int_{[0, 1]^d} \varv \cdot dDu .
  \]
  We call $\fracdiv\, \varv$ the divergence of $\varv$. This
  terminology is justified by the fact that
  \begin{equation}
    \label{eq:fracdivC}
  \forall u \in BV([0, 1]^d), \qquad (\fracdiv \, \varv)(u) =
  \int_{[0, 1]^d} u \operatorname{div} \varv
  \end{equation}
  whenever $\varv \in C^1([0, 1]^d; \R^d)$. To
  prove~\eqref{eq:fracdivC}, one may extend $\varv$ to a compactly
  supported $C^1$ vector field on $\R^d$, and then apply the trace
  theorem for a domain $U$ that is a bounded Lipschitz open
  neighborhood of $[0, 1]^d$ containing the support of $\varv$.

  Let us check that $\fracdiv \, \varv \in SCH([0, 1]^d)$ whenever
  $\varv \in C([0, 1]^d; \R^d)$. For $\varepsilon > 0$, choose a
  smooth vector field $w \colon \R^d \to \R^d$ such that $|\varv(x) -
  w(x)| \leq \varepsilon$ for all $x \in [0, 1]^d$. Then
  \[
  (\fracdiv \, \varv)(u) = - \int_{[0, 1]^d} u \operatorname{div} w -
  \int_{[0, 1]^d} (\varv - w) \cdot dDu
  \]
  Call $g$ the restriction of $\operatorname{div} w$ to $[0, 1]^d$. The
  preceding equality implies that
  \[
  \left| \fracdiv \, \varv(u) - T_g(u) \right| = \left| \int_{[0, 1]^d}
  (\varv-w) \cdot dDu \right| \leq \varepsilon \|u\|_{BV}
  \]
  As $u$ is arbitrary, this ensures that $\fracdiv \, \varv$ is a
  continuous linear functional and that $\|\fracdiv \, \varv - T_g\|
  \leq C\varepsilon$. This shows that $\fracdiv \, \varv$ belongs to
  the closure of $T(L^d([0, 1]^d))$, \textit{i.e} it is a strong charge functional.
  
  Thus, we have defined a linear map $\fracdiv \colon C([0, 1]^d;
  \R^d) \to SCH([0, 1]^d)$. On top of that, it is continuous, as can be seen
  from the inequality
  \[
  |(\fracdiv \, \varv)(u)| \leq \|\varv\|_\infty \|Du\|([0, 1]^d) \leq  \|\varv\|_\infty \|u\|_{BV}
  \]
  In fact, each strong charge is the divergence of a continuous vector
  field, as we show now.
\end{Empty}

\begin{Theorem}[Representation of strong charge functionals]
  \label{thm:fracdivonto}
  The operator $\fracdiv$ is onto.
\end{Theorem}

\begin{proof}
  The following proof is based on the results of \cite{DePaPfef} and
  \cite{BourBrez}, adapted to the present formalism.  First, we prove
  that the range of $\fracdiv$ is dense in $SCH([0, 1]^d)$. Let $g
  \colon [0, 1]^d \to \R$ be a smooth function. We will prove that
  $T_g = \fracdiv \, \varv$ for some continuous vector field
  $\varv$. As $C^\infty([0, 1]^d)$ is dense in $L^d([0, 1]^d)$, this
  will establish our claim. Let $w$ be a solution of the Poisson equation
  $\Delta w = g$ on $[0, 1]^d$ and set $\varv = \nabla w$. As $\varv$
  is of class $C^\infty$ on $[0, 1]^d$, we have
  \[
  T_g(u) = \int_{[0, 1]^d} ug = \int_{[0, 1]^d} u \operatorname{div}
  \varv = (\fracdiv\, \varv)(u)
  \]
  for all $u \in BV([0, 1]^d)$, which shows that $T_g = \fracdiv \, \varv$.

  We will conclude the proof by showing that the range of $\fracdiv$
  is closed. By~\cite[Theorem~3.1.21]{Megg}, it suffices to show that the range of the adjoint map $\fracdiv^*$
  is closed in $SCH([0, 1]^d)^*$.

  In the following diagram, we let $M([0, 1]^d; \R^d)$ be the Banach
  space of $\R^d$-valued Borel measures on $[0, 1]^d$ normed by total
  variation. It is isomorphic to the dual of $C([0, 1]^d;
  \R^d)$, by the Riesz representation theorem. The vertical arrows
  $\Upsilon$ and $\Phi$ are the obvious isomorphisms and we define $U
  = \Phi^{-1} \circ \fracdiv^* \circ \Upsilon$, so that the diagram is
  commutative.
  \[
  \begin{tikzcd}
    SCH([0, 1]^d)^* \arrow[rr, "\fracdiv^*"] && C([0, 1]^d; \R^d)^* \\
    BV([0, 1]^d) \arrow[u, "\Upsilon"] \arrow[rr, "U"]
    && M([0, 1]^d; \R^d) \arrow[u, swap, "\Phi"]
  \end{tikzcd}
  \]
  For all $u \in BV([0, 1]^d)$, we have $(\fracdiv^* \circ \Upsilon) (u) =
  \Upsilon(u) \circ \fracdiv$. Evaluating at $\varv \in C([0, 1]^d;
  \R^d)$ yields
  \[
  (\fracdiv^* \circ \Upsilon)(u)(\varv) = \Upsilon(u) (\fracdiv \,
  \varv) = (\fracdiv \, \varv)(u) = - \int_{[0, 1]^d} \varv \cdot dDu.
  \]
  On the other hand,
  \[
  (\Phi \circ U)(u)(\varv) = \int \varv \cdot dU(u).
  \]
  Since the preceding equalities hold for every $\varv$, we deduce
  that the measure $U(u)$ is the restriction of $-Du$ to $[0
    ,1]^d$. Therefore, its total variation is $\|U(u)\| = \|Du\|([0,
    1]^d) = \|Du\|(\R^d)$. We recall from
  Paragraph~\ref{e:defBVclosed} that $u \mapsto \|Du\|(\R^d)$ is a
  norm on $BV([0, 1]^d)$ that is equivalent to $\| \cdot \|_{BV}$.
  This completes the proof.
\end{proof}

\begin{Empty}[Remark in dimension 1]
  \label{emp:pi}
  In dimension $d = 1$, the space $SCH([0, 1])$ is isometric to
  $C_0([0, 1])$. Indeed, we already know that the continuous map
  $\fracdiv \colon C_0([0, 1]) \to SCH([0, 1])$ is onto. We define the map
  \[
  \Pi \colon SCH([0, 1]) \to C_0([0, 1])
  \]
  that sends $\alpha$ to $x \mapsto \alpha(\ind_{[0, x]})$. We claim
  that $\Pi \circ \fracdiv = \operatorname{id}$ and, thus, that
  $\fracdiv$ is an isomorphism whose inverse is $\Pi$. Indeed, letting
  $\varv \in C_0([0, 1])$ we have for all $x \in [0, 1]$,
  \[
  \Pi(\fracdiv \, \varv)(x) = (\fracdiv \, \varv)(\ind_{[0, x]}) = -
  \int_{[0, 1]} \varv \, dD\ind_{[0,x]} = \varv(x) - \varv(0) =
  \varv(x)
  \]
\end{Empty}
  
\begin{Empty}[Strong charges]
  \label{strong_charges}
  To each strong charge functional $\alpha$ we associate the map
  $\calS(\alpha) \colon \niceBV([0, 1]^d) \to \R$
  \[
  \calS(\alpha) : B \mapsto \alpha(\ind_B).
  \]
  It is clearly finitely additive, by the linearity of $\alpha$, and
  continuous with respect to $BV$-convergence, by
  Proposition~\ref{prop:weak*conv}. Thus, we have defined a linear map
  $\calS \colon SCH([0, 1]^d) \to CH([0, 1]^d)$. Charges that belong
  to the range of $\calS$ are called \emph{strong}, for historical
  reasons. We admit for now that $\calS$ is injective -- a fact that
  will be proved momentarily in Corollary~\ref{cor:Sinjective}. As a consequence, strong
  charge functionals and strong charges are basically the same
  objects, under different disguises. However, the formalism of strong
  charge functionals makes it readily clear that $SCH([0, 1]^d)$ is a
  Banach space and this allows us to introduce a functional
  analytic approach in the next section.

  In dimension $d = 1$, charges on $[0, 1]$ are automatically
  strong. Indeed, let $\mu$ be any charge on $[0, 1]$ and $\varv
  \colon x \mapsto \mu([0, x])$ its associated continuous function,
  then it is easy to check that $\mu = \calS(\fracdiv \, \varv)$.
  
  The structure of $BV$-sets in $[0, 1]^d$ is fully understood. They
  have a Borel measurable {\it reduced} boundary $\partial_* B$, defined in
  a measure-theoretic way, on which a normal unit outer vector field
  $n_B$ is defined $\calH^{d-1}$-almost everywhere and they satisfy a
  generalized Gauss-Green theorem~\cite[5.9, Theorem~1]{EvanGari}
  \[
  \int_B \operatorname{div} \varphi = \int_{\partial_*B} \varphi \cdot
  n_B \, d\calH^{d-1},
  \]
  for all $\varphi \in C^1([0, 1]^d; \R^d)$. In particular, if a
  strong charge functional $\alpha$ is the divergence of $\varv \in
  C([0, 1]^d; \R^d)$, a density argument shows that
  \[
  \alpha(\ind_B) = \int_{\partial_* B} \varv \cdot n_B \,
  d\calH^{d-1}.
  \]
  In terms of charges, a charge $\mu$ is strong if and only if there exists a continuous vector field $\varv \in C([0, 1]^d; \R^d)$ such that
  \[
  \mu(B) = \int_{\partial_* B} \varv \cdot n_B \, d\calH^{d-1}, \qquad B \in \niceBV([0, 1]^d).
  \]
  A function $f \in C_0([0, 1]^d)$ is called \emph{strongly
    chargeable} whenever $\Delta_f$ extends (uniquely) to a strong
  charge. For example, any continuous function $f \in C_0([0, 1])$ is
  strongly chargeable. This only happens in dimension $d=1$.
\end{Empty}

We close this section with an approximation lemma whose full strength
will prove useful in Section~\ref{sec:Faber}, in the proof of
Theorem~\ref{thm:mainSchauder}. A dyadic partition $\calP$ of $[0,
  1]^d$ is a finite set of pairwise almost disjoint dyadic cubes in
$I=[0, 1]^d$ such that $\bigcup \{ K : K \in \calP \} = [0, 1]^d$. We
do not require that the cubes in $\calP$ are of the same generation.

\begin{Lemma}
  \label{thm:approx2}
  Let $(\calP_n)$ be a sequence of dyadic partitions of $[0, 1]^d$ whose
  meshes tend to~$0$, i.e
  \[
  \lim_{n \to \infty} \max \left\{ \rmdiam K : K \in \calP_n \right\} = 0.
  \]
  Let $u \in BV([0, 1]^d)$.  For each $n$, define the
  function
  \[
  u_n = \sum_{K \in \calP_n} \bar{u}_K \ind_{K}, \qquad \text{where }
  \bar{u}_K = \frac{1}{|K|} \left(\int_K u \right).
  \]
  Then $u_n \to u$ weakly* in $BV([0, 1]^d)$.
\end{Lemma}

\begin{proof}
  In regard with Proposition~\ref{prop:weak*conv}, we need to prove
  that $u_n \to u$ in $L^1$ and that $\sup_n \|Du_n\|(\R^d) <
  \infty$. That $u_n \to u$ in $L^1$ is routinely proven by
  approximating $u$ by a continuous function $\varv$ on $[0, 1]^d$ and
  using the uniform continuity of $\varv$.
  
  Thus, we concentrate our efforts on the second part. By the
  continuity of the trace operator, there is a constant $C \geq 0$ such that
  \[
  \int_{\partial ((0, 1)^d)} |\rmtr \varphi| d\calH^{d-1} \leq C \left(\int_{(0, 1)^d}
  |\varphi| + \|D\varphi\|((0, 1)^d) \right),
  \]
  for all $\varphi \in BV((0, 1)^d)$. Then a scaling argument shows that
  \[
  \int_{\partial K} |\rmtr \varphi| d\calH^{d-1} \leq \frac{C}{\operatorname{diam} K} 
  \int_{\operatorname{int} K} |\varphi| + C \|D\varphi \|(\operatorname{int} K)
  \]
  whenever $K \subset [0, 1]^d$ is a dyadic cube and $\varphi \in BV(\operatorname{int} K)$.

  Now we fix $n$. For each cube $K \in \calP_n$, we call
  $u_{\mid\operatorname{int}K} \in BV(\operatorname{int} K)$ the
  restriction of $u$ to $\operatorname{int} K$ and $\varv_K = (u -
  \bar{u}_K) \ind_K \in BV([0, 1]^d)$.  By the extension theorem,
  \[
  \|D\varv_K\|(\R^d) = \|Du\|(\operatorname{int} K) + \int_{\partial
    K} |\rmtr (u_{\mid\operatorname{int}K} - \bar{u}_K)| \, d\calH^{d-1}.
  \]
  Hence, with the help of the preceding inequality and the Hölder
  inequality, we deduce that
  \begin{align*}
    \|D\varv_K\|(\R^d) & \leq (1+ C) \|Du\|(\operatorname{int} K) +
    \frac{C}{\operatorname{diam} K} \int_{K} |u - \bar{u}_K| \\ &
    \leq (1 + C) \|Du\|(\operatorname{int} K) + \frac{C}{\sqrt{d}}
    \left(\int_{K} |u - \bar{u}_K|^{d/(d-1)} \right)^{1 - 1/d}.
  \end{align*}
  By the Sobolev-Poincaré inequality on a cube, one has
  \[
  \left(\int_{K} |u - \bar{u}_K|^{d/(d-1)} \right)^{1 - 1/d} \leq C_P
  \|Du\|(\operatorname{int} K).
  \]
  where $C_P$ is the Poincaré constant, hence,
  \begin{equation}
    \label{eq:jetrouvepasdenom}
  \|D\varv_K\|(\R^d) \leq C' \|Du\|(\operatorname{int} K)
  \end{equation}
  for some constant $C'$ (depending solely on $d$). Finally, one
  notices that 
  \begin{align*}
  \|D u_n\|(\R^d) & \leq \left\| D \sum_{K \in \calP_n} \varv_K
  \right\|(\R^d) + \|Du\|(\R^d) \\ & \leq \sum_{K \in \calP_n}
  \|D\varv_K\|(\R^d) + \|Du\|(\R^d) \\ & \leq (1 +C') \|Du\|(\R^d)
  \end{align*}
  In the last inequality, we used~\eqref{eq:jetrouvepasdenom} and the
  fact that the interiors of the cubes $K \in \calP_n$ are pairwise
  disjoint.
\end{proof}

\begin{Corollary}
  \label{cor:Sinjective}
  The map $\calS \colon SCH([0, 1]^d) \to CH([0, 1]^d)$ is injective.
\end{Corollary}

\begin{proof}
  Let $\alpha \in SCH([0, 1]^d)$ be in the kernel of $\calS$ and $u
  \in BV([0, 1]^d)$. Then, for every dyadic cube $K$, one has
  $\alpha(\ind_K) = \calS(\alpha)(K) = 0$. For any integer $n$,
  consider the collection $\calP_n$ of all dyadic cubes of generation
  $n$ and define the sequence of approximating functions $u_n =
  \sum_{K \in \calP_n} \bar{u}_K \ind_K$ as in the preceding
  lemma. Then $\alpha(u_n) = \sum_{K\in \calP_n} \bar{u}_K
  \alpha(\ind_K) = 0$, for all $n$. Moreover, $u_n \to u$ in the weak*
  topology. But $\alpha$ belongs to the predual of $BV([0, 1]^d)$,
  consequently, $\alpha(u) = \lim \alpha(u_n) = 0$. This is true for
  any $u$, so, we infer that the kernel of $\calS$ is trivial.
\end{proof}

\section{The Faber-Schauder basis in $SCH([0, 1]^d)$}
\label{sec:Faber}

\begin{Empty}[Schauder basis]
  We recall that a sequence $(e_n)$ with terms in an
  infinite-dimensional Banach space $X$ is called a Schauder basis
  (or simply basis) of $X$ if for each $x \in X$, there is a unique
  sequence $(a_n)$ of scalars such that
  \begin{equation}
    \label{eq:Schauder}
    x = \sum_{n = 0}^\infty a_n e_n
  \end{equation}
  strongly in $X$.

  If the convergence in~\eqref{eq:Schauder} is required to be weak
  instead of strong, then we say that $(e_n)$ is a weak basis. A
  result of Mazur states that a weak basis is in fact a Schauder
  basis, see \cite[Theorem 5.3]{Morr} where this result is attributed
  to Banach.
\end{Empty}

\begin{Empty}[Haar basis]
  An example of a Schauder basis is provided by the system of Haar
  functions described here. It is a basis in every $L^p([0, 1])$ space,
  where $1 \leq p < \infty$. Set the functions $h_{-1} = 1$,
  \[
  h_{0, 0} \colon x \in [0, 1] \mapsto \begin{cases}
    1 & \text{if } 0 \leq x < 1/2 \\
    -1 & \text{if } 1/2 \leq x \leq 1
  \end{cases}
  \]
  (one may extend $h_{0, 0}$ to $\R$ by zero if needed in the
  subsequent formulae) and then, for every integer $n \geq 1$ and $k =
  0, \dots, 2^{n} - 1$,
  \[
  h_{n,k} \colon x \mapsto 2^{n/2} h_{0, 0}(2^{n} x - k) = \begin{cases}
    2^{n/2} & \text{if } 2^{-n} k \leq x < 2^{-n}k + 2^{-n-1} \\
    -2^{n/2} & \text{if } 2^{-n}k + 2^{-n-1} \leq x < 2^{-n}(k+1) \\
    0 & \text{otherwise}.
  \end{cases}
  \]
  The Haar basis is the sequence $h_{-1}, h_{0, 0}, h_{1,0}, h_{1, 1},
  h_{2, 0},h_{2,1},h_{2,2}, h_{2,3}, \dots$ (indices are ordered lexicographically). With our
  normalization choice this is an orthonormal basis in $L^2([0, 1])$.
\end{Empty}

\begin{Empty}[One-dimensional Faber-Schauder system]
  \label{emp:faber1d}
  By definition, the Faber-Schauder functions are the primitives of
  the Haar functions, that is,
  \[
  f_{-1} \colon x \mapsto \int_0^x h_{-1}, \qquad f_{n,k} \colon x \mapsto
  \int_0^x h_{n,k}
  \]
  for $n \geq 0$ and $k = 0, \dots, 2^{n}-1$. It was first proven in
  \cite{Schau} that $f_{-1}, f_{0,0}, f_{1,0}, f_{1,1}, f_{2,0},
  \dots$ constitutes a Schauder basis of $C_0([0, 1])$.  We recall our
  claim that the map $\Pi$ introduced in~\ref{emp:pi} is an
  isomorphism between $SCH([0, 1])$ and $C_0([0, 1])$.  As $f_{-1} =
  \Pi(T_{h_{-1}})$ and $f_{n,k} = \Pi(T_{h_{n,k}})$ for all indices
  $n, k$, we can assert that $T_{h_{-1}}$, $T_{h_{0,0}}$,
  $T_{h_{1,0}}$, $T_{h_{1,1}}$, $T_{h_{2,0}}, \dots$ is a Schauder
  basis of $SCH([0, 1])$.
\end{Empty}

\begin{Empty}[Multidimensional Haar basis]
  \label{emp:multiHaar}
  Define the matrix
  \[
  A = \begin{pmatrix}
    1 & 1 \\ 1 & -1
  \end{pmatrix}
  \]
  and let $A_d$ be the matrix of order $2^d$ that is the
  Kronecker product of $A$ with itself $d$ times, \textit{i.e.}
  $A_1 = A$ and
  \[
  A_{d+1} = \begin{pmatrix} A_{d} & A_d \\ A_d & -A_d
  \end{pmatrix} ,\qquad \text{for } d \geq 1 .
  \]
  By induction, $A_d$ is easily seen to be a symmetric matrix such
  that $(A_d)^2 = 2^d I_{2^d}$ (where $I_{2^d}$ denotes the identity
  matrix of order $2^d$). In other words, $2^{-d/2}A_d$ is an orthogonal matrix.
  Subsequently, the entries of the matrix
  $A_d$ will be written $(A_d)_{r, \ell}$, where the row and column
  numbers $r$ and $\ell$ range over $\{0, \dots, 2^d - 1\}$.

  For all $n \geq 0$, we let $K_{n, k}$, $k = 0, \dots, 2^{nd} - 1$, be
  the collection of all dyadic cubes in $[0, 1]^d$ of side-length
  $2^{-n}$. We further require that, for all indices $n$ and $k$, the
  cubes $K_{n+1, 2^d k}$, $K_{n+1, 2^d k + 1}, \dots K_{n+1, 2^dk +
    2^d -1}$ are the $2^d$ subcubes of $K_{n, k}$ of side
  $2^{-(n+1)}$.

  We are now able to build the Haar basis. The first Haar function is
  the exceptional one
  \[
  g_{-1} \colon (x, y) \in [0, 1]^d \mapsto 1.
  \]
  Then, we define, for all $n \geq 0$, $k \in \{0, \dots, 2^{nd} -
  1\}$ and $r \in \{1, \dots, 2^d - 1\}$:
  \begin{equation}
    \label{eq:defgnkr}
    g_{n, k, r} = 2^{nd/2} \sum_{\ell=0}^{2^d - 1} (A_d)_{r, \ell}
    \ind_{K_{n+1, 2^d k + \ell}}.
  \end{equation}
  We refer to $n$, $k$ and $r$ as the generation number, the cube
  number and the type number of $g_{n, k, r}$. By construction,
  $g_{n,k,r} = 0$ almost everywhere outside the cube $K_{n,k}$. In fact,
  the support of $g_{n,k,r}$ is $K_{n,k}$ and the support of $g_{-1}$ is $K_{0,0}$. We
  also note that the average value of $g_{n,k,r}$ is $0$ (this is
  because the $r$-th line of $A_d$ is orthogonal to the zeroth line,
  which is filled with ones).

  The Haar basis is $g_{-1}, g_{0, 0, 1}, \dots, g_{0,0, 2^{d}-1},
  g_{1, 0, 1}, \dots, g_{1, 2^d-1, 2^d -1}, g_{2, 0, 1}, \dots$
  (indices are ordered lexicographically).
  \begin{Claim*}
    The Haar functions are orthonormal in $L^2([0, 1]^d)$.
  \end{Claim*}
  \noindent Those functions are indeed appropriately normalized; and
  the cases worth considering in proving that these functions are
  pairwise orthogonal are:
  \begin{itemize}
  \item The case of two functions $g_{n, k, r}$ and $g_{n', k', r'}$
    with $n < n'$: if $K_{n',k'} \subset K_{n, k}$, then
    \[
    \int g_{n,k,r} g_{n', k', r'} = \pm \int g_{n', k', r'} = 0.
    \]
    Otherwise, $|K_{n', k'} \cap K_{n,k}| = 0$, thus, $g_{n,k,r}
    g_{n', k', r'} = 0$ a.e., which implies that $g_{n,k,r}$ and $g_{n',
      k', r'}$ are orthogonal.
  \item The case of two functions $g_{n, k, r}$ and $g_{n, k, r'}$
    of the same generation and same cube numbers but different type
    numbers $r \neq r'$:
    \[
    \int g_{n,k,r} g_{n', k', r'} = \frac{1}{2^d} \sum_{\ell = 0}^{2^d
      -1} (A_d)_{r, \ell} (A_d)_{r', \ell} = 0.
    \]
    Here, we use that the matrix $2^{-d/2} A_d$ is orthogonal.
  \end{itemize}
  \begin{Claim*}
    The function $\ind_{K_{n,k}}$ is a linear combination of $g_{-1}$
    and the functions $g_{n', k', r}$ corresponding to generation numbers 
    $n' \in \{0,1,\ldots,n-1\}$,
    cube numbers $k' = \left\lfloor 2^{-d(n-n')} k \right\rfloor$ (this is equivalent to
    $K_{n',k'} \supset K_{n,k}$), and
    arbitrary type numbers $r \in \{1, \dots, 2^d - 1\}$.
  \end{Claim*}
  \noindent This claim is proven by induction on $n$. The base case is
  straightforward as $\ind_{K_{0,0}} = g_{-1}$. Regarding the induction
  step, we note that
  \begin{equation}
    \label{eq:relationCubesHaar}
  \begin{pmatrix}
    2^{nd/2} \ind_{K_{n,k}} \\ g_{n, k, 1} \\ \vdots \\ g_{n, k, 2^d
      -1}
  \end{pmatrix} =
  2^{nd/2} A_d \begin{pmatrix} \ind_{K_{n+1, 2^d k}} \\ \ind_{K_{n+1,
        2^d k + 1}} \\ \vdots \\ \ind_{K_{n+1, 2^d k + 2^d - 1}}
  \end{pmatrix}
  \end{equation}
  The conclusion follows from the invertibility of $A_d$.
  
  The main result of this section is that, by primitivizing (i.e., by
  applying $T$ to) the Haar
  functions in the space of strong charge functionals, we obtain a Schauder basis
  of $SCH([0, 1]^d)$, see Theorem~\ref{thm:mainSchauder}. In view
  of~\ref{emp:faber1d}, this basis is analogous to the 1-dimensional
  Faber-Schauder basis. It may be worth noting that, like bases of wavelets, the
  Faber-Schauder basis we obtain is ``localized'' in the sense that the supports of
  its members are controlled; in fact, $\mathrm{supp}\, T_{g_{n,k,r}} = K_{n,k}$.

  We warn the reader that a Schauder basis does not need to be unconditional,
  that is, the order of summation matters. In fact, the
  $1$-dimensional Faber-Schauder basis of $C_0([0, 1])$ is not
  unconditional, see \cite{AlbiKalt}.
  %First, we need the following proposition, of independent
  %interest. It intuitively relates to the fact that the charge
  %functionals $T_{g_{n,k}}$, $k = 0, \dots, 3\cdot 4^n - 1$ have
  %``localized supports''.
\end{Empty}

\begin{Theorem}
  \label{thm:mainSchauder}
  The sequence $T_{g_{-1}}, T_{g_{0, 0, 1}}, T_{g_{0, 0, 2}}, \dots$
  is a Schauder basis of the space $SCH([0, 1]^d)$, with respect to which each strong charge functional
  $\alpha$ is decomposed as follows:
  \[
  \alpha = \alpha(g_{-1}) T_{g_{-1}} + \alpha(g_{0, 0, 1}) T_{g_{0, 0,
      1}} + \alpha(g_{0, 0, 2}) T_{g_{0, 0, 2}} + \cdots .
  \]
  (The convergence occurs, of course, in $BV([0,1]^d)^*$.)
\end{Theorem}

\begin{proof}
  First, we prove the uniqueness part in the definition of a Schauder
  basis. Suppose a strong charge functional $\alpha$ has a decomposition
  \begin{multline}
    \label{eq:decompSchauder}
    \alpha = a_{-1} T_{g_{-1}} + a_{0, 0, 1} T_{g_{0, 0, 1}} + a_{0,
      0, 2} T_{g_{0, 0, 2}} + \cdots + a_{0, 0, 2^{d}-1} T_{g_{0, 0,
        2^{d} - 1}} \\ + a_{1, 0, 1} T_{g_{1, 0, 1}} + \cdots + a_{1,
      2^{d}-1, 2^{d}-1} T_{g_{1, 2^{d}-1, 2^{d} -1}} + \cdots
  \end{multline}
  Applying~\eqref{eq:decompSchauder} at $g_{n,k,r}$, we get $a_{n,k,r}
  = \alpha(g_{n,k,r})$ by the orthonormality of the Haar
  functions. Likewise, $a_{-1} = \alpha(g_{-1})$.
  
  We turn to the existence. We fix a strong charge functional $\alpha \in
  SCH([0, 1]^d)$. Our goal is to prove that~\eqref{eq:decompSchauder} holds
  weakly, when the coefficients $a_{-1}$ and $a_{n,k,r}$ are taken as
  in the first part of the proof. Mazur's weak basis theorem will then
  imply the desired result. Therefore, we consider a function $u \in
  BV([0, 1]^d)$ and wish to prove that
  \begin{multline}
    \label{eq:decompSchauder2}
    \alpha(u) = a_{-1} T_{g_{-1}}(u) + a_{0, 0, 1} T_{g_{0, 0, 1}}(u)
    + a_{0, 0, 2} T_{g_{0, 0, 2}}(u) + \cdots + a_{0, 0, 2^{d}-1}
    T_{g_{0, 0, 2^{d} - 1}}(u) \\ + a_{1, 0, 1} T_{g_{1, 0, 1}}(u) +
    \cdots + a_{1, 2^{d}-1, 2^{d}-1} T_{g_{1, 2^{d}-1, 2^{d} -1}}(u) +
    \cdots .
  \end{multline}
  For all $n, k$, we define the truncated sum 
  \[
  \alpha_{n,k} = a_{-1} T_{g_{-1}} + \cdots + a_{n, k, 2^d - 1}
  T_{g_{n, k, 2^d -1}}.
  \]
  We will prove that $\alpha_{n, k}(u) \to \alpha(u)$ (as usual, we equip the set of
  legal couples of indices $(n,k)$ with the lexicographical order). The sequence
  $(\alpha_{n,k}(u))$ is merely a subsequence of the sequence of partial
  sums in~\eqref{eq:decompSchauder2} and we shall deal with this issue
  at the end of this proof. Define
  \[
  \calG_{n,k} = \left\{ g_{-1}, g_{0, 0, 1}, \dots, g_{n-1,2^{(n-1)d},2^d-1},g_{n,0,1},
  		\ldots,g_{n, k,
    2^d - 1} \right\}
  \]
  so that
  \[
  	\alpha_{n,k} = \sum_{g \in
    \calG_{n,k}} \alpha(g) T_g.
  \]
  Define also the dyadic partition of $[0,1]^d$
  \[
  \calP_{n,k} = \left\{K_{n, k'} : k+1 \leq k' \leq 2^{nd} - 1\right\}
  \cup
  \left\{ K_{n+1, k'} : 0 \leq k' \leq 2^d k + 2^d - 1
  \right\}.
  \]
  Reasoning as in the first part of the proof, $\alpha_{n,k}(g) =
  \alpha(g)$ for all $g \in \calG_{n,k}$. Also, with the help of the
  second claim of~\ref{emp:multiHaar}, we have
  \begin{equation}
    \label{eq:spanGnk}
  \operatorname{span} \calG_{n,k} \supset \operatorname{span} \left\{ \ind_K
  : K \in \calP_{n,k} \right\}.
  \end{equation}
  The two points above guarantee that $\alpha(\ind_K) =
  \alpha_{n,k}(\ind_K)$ for all $K \in \calP_{n,k}$.

  Next, we define
  \[
  u_{n,k} = \sum_{K \in \calP_{n,k}} \frac{1}{|K|} \left(\int_K
  u\right) \ind_{K}.
  \]
  As each function $g \in \calG_{n,k}$ is constant a.e. on the dyadic cubes
  in $\calP_{n,k}$, we clearly
  have $T_g(u) = T_g(u_{n,k})$. Therefore,
  \begin{align*}
    \alpha_{n,k}(u) & = \alpha_{n,k}(u_{n,k}) \\ & = \sum_{K \in
      \calP_{n,k}} \frac{1}{|K|} \left( \int_K u \right)
    \alpha_{n,k}(\ind_K) \\ & = \sum_{K \in \calP_{n,k}} \frac{1}{|K|}
    \left( \int_K u \right) \alpha(\ind_K) \\ & = \alpha(u_{n,k}).
  \end{align*}
  Lemma~\ref{thm:approx2} applies, therefore, $u_{n,k} \to u$ weakly*
  by Proposition~\ref{prop:weak*conv}, from which we deduce that
  $\alpha(u_{n,k}) \to \alpha(u)$.
  
  To finish the proof, we ought to show that the sequence of partial sums in the
  right-hand side of~\eqref{eq:decompSchauder2} tends to $0$. Any partial sum in
  this sequence differs from some $\alpha_{n,k}(u)$ considered above by at most
  $2^d -1$ terms of the type 
  \[
  	a_{n,k,r}T_{g_{n,k,r}}(u) = \alpha(g_{n,k,r})\int g_{n,k,r} u = 
  	\alpha \left( \left( \int g_{n,k,r}u\right) g_{n,k,r} \right).
  \] 
  In other words, we must establish that for all $\alpha \in SCH([0,1]^d)$ and
  all $u \in BV([0,1]^d)$, $a_{n,k,r}T_{g_{n,k,r}}(u) \to 0$. Fixing $u \in BV([0,1]^d)$,
  this is equivalent to showing that the sequence $\left( \left( \int g_{n,k,r}u
  \right) g_{n,k,r} \right)$ weakly* converges to 0 (with respect to the duality set forth 
  in Theorem \ref{thm:predual}), according to the equation above. By proposition
  \ref{prop:weak*conv}, this is in turn equivalent to showing that 
  \begin{equation}
  \label{eq:th-1}
  \left\|  \left( \int g_{n,k,r}u \right) g_{n,k,r}\right\|_1 \to 0
  \end{equation}
  and
  \begin{equation}
  \label{eq:th-2}
  \sup \left\|  \left( \int g_{n,k,r}u \right) g_{n,k,r}\right\|_{BV} < \infty.
  \end{equation}
  When establishing these two facts we will distinguish between the case $d=1$ and
  the case $d \geq 2$. It will be useful to note that: for all $d$, $\|g_{n,k,r}\|_1 = 
  2^{nd/2}|K_{n,k}| = 2^{nd/2}2^{-nd} = 2^{-nd/2}$; if $d=1$ then $\|u\|_\infty < \infty$;
  if $d \geq 2$ then $\|g_{n,k,r}\|_{d/(d-1)} = 2^{nd/2}|K_{n,k}|^{(d-1)/d} = 
  2^{nd/2}2^{-n(d-1)} = 2^{n(1-d/2)}$.
  
  Thus, in case $d=1$, one has
  \[
  \left\|  \left( \int g_{n,k,r}u \right) g_{n,k,r}\right\|_1 \leq 
  \left( \int |g_{n,k,r}| |u| \right) \|g_{n,k,r}\|_1 \leq  
  \|u\|_\infty \|g_{n,k,r}\|_1^2 \leq
  \|u\|_\infty 2^{-n} \to 0,
  \]
  whereas in case $d \geq 2$ one has
  \begin{multline*}
  \left\|  \left( \int g_{n,k,r}u \right) g_{n,k,r}\right\|_1 \leq 
  \left( \int |g_{n,k,r}| |u| \right) \|g_{n,k,r}\|_1 \leq  
  \|u\|_d \|g_{n,k,r}\|_{d/(d-1)} \|g_{n,k,r}\|_1 \\
  = \|u\|_d 2^{n(1-d/2)} 2^{-nd/2} = \|u\|_d 2^{n(1-d)}\to 0.
  \end{multline*}
  Accordingly, \eqref{eq:th-1} holds for all $d$.
  
  Next, it is useful to note that: if $d=1$ then $\|g_{n,k,r}\|_{BV} \leq 2.2^{n/2}.4$
  (the height of each jump is $2.2^{n/2}$ and there are maximum 4 jumps); if $d\geq 2$ then
  $ \|g_{n,k,r}\|_{BV} \leq 2^{nd/2}2^d2(2d)2^{-(n+1)(d-1)} = 8d2^{n(1-d/2)}$.
  
  Thus, in case $d=1$ one has
  \begin{multline*}
  \left\|  \left( \int g_{n,k,r}u \right) g_{n,k,r}\right\|_{BV} \leq 
  \left( \int |g_{n,k,r}| |u| \right) \|g_{n,k,r}\|_{BV} \leq  
  \|u\|_\infty \|g_{n,k,r}\|_1 \|g_{n,k,r}\|_{BV} \\
  \leq \|u\|_\infty 2^{-n/2}8.2^{n/2} \leq 8 \|u\|_\infty,
  \end{multline*}
  whereas in case $d \geq 2$ one has
  \begin{multline*}
  \left\|  \left( \int g_{n,k,r}u \right) g_{n,k,r}\right\|_{BV} \leq 
  \left( \int |g_{n,k,r}| |u| \right) \|g_{n,k,r}\|_{BV} \\  
  \leq \|u\|_d \|g_{n,k,r}\|_{d/(d-1)} \|g_{n,k,r}\|_{BV} 
  = \|u\|_d 2^{n(1-d/2)}8d2^{n(1-d/2)} \\
  = 8d\|u\|_d 2^{n(2-d)} \leq 8d\|u\|_d.
  \end{multline*}
  Accordingly, \eqref{eq:th-2} holds for all $d$. This completes the proof.
\end{proof}

\begin{Empty}[Remark on charge functionals]
  \label{e:FaberCH}
  We briefly outline how it is possible to define a notion of charge
  functional similar to that of strong charge functional, thereby
  endowing $CH(A)$ with a Banach space structure. First, define the
  space $BV_\infty(A)$ of measurable functions $u \colon \R^d \to \R$
  that are essentially bounded, that vanish almost everywhere outside of $A$, 
  and that have bounded varition. This space is normed
  by $\|u\|_{BV_\infty} = \|u\|_\infty +
  \|Du\|(\R^d)$. We define the linear map $T \colon L^1(A) \to
  BV_\infty(A)^*: f \mapsto T_f$ by
  \[
  T_f(u) = \int_A fu, \qquad f \in L^1(A), \;u \in BV_\infty(A).
  \]
  The space of charge functionals is the closure of $T(L^1(A))$ in
  $BV_\infty(A)^*$. This space is a canonical predual of
  $BV_\infty(A)$. Charge functionals are in bijection
  with charges on $A$: if $\alpha$ is a charge functional, then $B
  \mapsto \alpha(\ind_B)$ is a charge. In case $A = [0, 1]^d$, the
  image of the Haar basis of $L^1([0, 1]^d)$ under $T$ is a Schauder
  basis of the space of charge functionals, similar to the
  Faber-Schauder basis of $SCH([0, 1]^d)$.
\end{Empty}

\section{Criteria for strong chargeability}

\begin{Empty}
  Let $f \in C_0([0, 1]^d)$, we define $\lambda_{-1}(f) = \Delta_f([0,
    1]^d) = f(1, \dots, 1)$ and, for all relevant indices $n, k$ and
  $r$,
  \begin{equation}
    \label{eq:formulaeLambda}
    \lambda_{n,k,r}(f) = 2^{nd/2} \sum_{\ell=0}^{2^d - 1} (A_d)_{r,
      \ell} \Delta_f (K_{n+1, 2^d k + \ell}) 
  \end{equation}
  The maps $\lambda_{-1}$ and $\lambda_{n, k, r}$ defined above are
  continuous linear functionals on $C_0([0, 1]^d)$. It is clear that
  if $f$ is strongly chargeable, then $\Delta_f$ is by definition a
  strong charge and $\lambda_{n, k, r}(f)$ are the coefficients of the
  strong charge functional $\calS^{-1}(\Delta_f)$ in the
  Faber-Schauder basis. From this observation we derive the
  Theorem~\ref{thm:chargeability}, that equates strong chargeability
  with the convergence of a series in $SCH([0, 1]^d)$.
\end{Empty}

\begin{Theorem}
  \label{thm:chargeability}
  A function $f \in C_0([0, 1]^d)$ is strongly chargeable if and only
  if the Faber-Schauder series
  \begin{equation}
    \label{eq:decompChargeability}
    \lambda_{-1}(f) T_{g_{-1}} + \lambda_{0, 0, 1}(f) T_{g_{0, 0, 1}}
    + \cdots
  \end{equation}
  converges in $SCH([0, 1]^d)$.
\end{Theorem}

\begin{proof}
  The direct implication follows from the arguments above. Conversely,
  suppose the series~\eqref{eq:decompChargeability} is convergent and
  denote its sum by $\alpha$. Then
  \[
  \alpha(g_{-1}) = \alpha(\ind_{K_{0,0}}) = \lambda_{-1}(f) =
  \Delta_f (K_{0,0}) \text{ and } \alpha(g_{n,k,r}) =
  \lambda_{n,k,r}(f)
  \]
  by Theorem~\ref{thm:mainSchauder}.

  We now prove by induction on $n$ that $\alpha(\ind_{K_{n,k}}) =
  \Delta_f (K_{n,k})$ for all $k = 0, \dots, 2^{nd} - 1$. The base case
  $n = 0$ is already treated. Suppose the result is valid for a
  generation number $n \geq 0$. Fix a cube number $k$. By applying
  $\alpha$ to~\eqref{eq:relationCubesHaar}, we get
  \[
  \begin{pmatrix}
    2^{nd/2} \Delta_f (K_{n,k}) \\
    \lambda_{n,k,1}(f) \\
    \vdots \\
    \lambda_{n,k, 2^d - 1}(f)
  \end{pmatrix} = 2^{nd/2} A_d
  \begin{pmatrix}
    \alpha(\ind_{K_{n+1, 2^d k}}) \\
    \alpha(\ind_{K_{n+1, 2^d k + 1}}) \\
    \vdots \\
    \alpha(\ind_{K_{n+1, 2^d k + 2^d - 1}})
  \end{pmatrix}
  \]
  On the other hand, from the definition~\eqref{eq:formulaeLambda} of
  the $\lambda_{n,k,r}$ functionals, we have
  \[
  \begin{pmatrix}
    2^{nd/2} \Delta_f (K_{n,k}) \\
    \lambda_{n,k,1}(f) \\
    \vdots \\
    \lambda_{n,k, 2^d - 1}(f)
  \end{pmatrix}
  = 2^{nd/2} A_d
  \begin{pmatrix}
    \Delta_f (K_{n+1, 2^d k}) \\ \Delta_f (K_{n+1, 2^d k + 1}) \\ \vdots
    \\ \Delta_f (K_{n+1, 2^d k + 2^d - 1})
  \end{pmatrix}.
  \]
  As $A_d$ is invertible, this ends the proof by induction.

  Let $(x_1, \dots, x_d) \in [0, 1]^d$ a point whose coordinates are
  dyadic numbers. Then we can write $K = \prod_{i=1}^d [0, x_i]$ as a
  finite union of almost disjoint dyadic cubes. Using the result
  above, we derive that $\calS^{-1}(\alpha)(K) = \alpha(\ind_K) =
  \Delta_f (K) = f(x_1, \dots, x_d)$.  When $(x_1, \dots, x_d) \in [0,
    1]^d$ is arbitrary, we use a simple density argument (and the
  continuity of $f$ and $\alpha$) to justify that
  $\calS^{-1}(\alpha)(K) = f(x_1, \dots, x_d)$ holds as well. Therefore,
  $f$ is chargeable, by Proposition~\ref{prop:primcharge}, and $\Delta_f
  = \calS^{-1}(\alpha)$, which ensures that $f$ is strongly
  chargeable.
\end{proof}

\begin{Empty}
  The 1-dimensional Faber-Schauder functions $f_{n,k}$
  (see~\ref{emp:faber1d}) have the pleasant property of having
  localized supports. This helps estimate the norm of a linear
  combination of $f_{n,k}$ functions, for a fixed generation number $n
  \geq 0$. Indeed, we clearly have
  \[
  \left\| \sum_{k=0}^{2^n - 1} a_k f_{n, k} \right\|_\infty =
  \frac{1}{2^{n/2+1}}\max_{0 \leq k \leq 2^n - 1} |a_k|.
  \]
  The following Proposition~\ref{prop:estimate} is a subtler
  multidimensional analogue. Together with
  Theorem~\ref{thm:chargeability}, it will allow to state strong
  chargeability or non strong chargeability criteria of practical use
  in Corollary~\ref{cor:chargeability}.
\end{Empty}

\begin{Proposition}
  \label{prop:estimate}
  There is a positive constant $C$ such that for all $n \geq 0$ and
  scalars $(a_{k, r})$, we have
  \[
  \frac{1}{2^{n(d/2+1)} C} \max_{1 \leq r \leq 2^d -1} 
  \sum_{k=0}^{2^{nd}-1} |a_{k, r}| \leq \left\| \sum_{k=0}^{2^{nd} - 1}
  \sum_{r=1}^{2^d - 1} a_{k , r} T_{g_{n,k, r}} \right\| \leq C
  2^{n(d/2 - 1)} \max_{k, r} |a_{k, r}|
  \]
\end{Proposition}

\begin{proof}
  First we prove the upper bound. Let $u \in BV([0, 1]^d)$. For any cube
  number $k$ and type number $r$, we have
  \[
  a_{k,r} T_{g_{n, k, r}}(u) = a_{k,r} \int g_{n, k, r} u = a_{k, r}
  \int g_{n, k, r} (u - \bar{u}_{K_{n, k}})
  \]
  where $\bar{u}_{K_{n,k}}$ is the average value of $u$ on the $K_{n,
    k}$. By the Hölder inequality and the Poincaré inequality,
  we obtain that
  \begin{align*}
  \left|a_{k, r} T_{g_{n, k, r}}(u) \right| & \leq |a_{k, r}|
  \|g_{n,k,r}\|_d \left(\int_{K_{n,k}} |u - \bar{u}_{K_{n,k}}|^{d/(d-1)}
  \right)^{1 - 1/d} \\ & \leq C_P 2^{n(d/2 - 1)} |a_{k,r}| \,
  \|Du\|(\operatorname{int} K_{n,k})
  \end{align*}
  where $C_P$ is the Poincaré constant for $d$-dimensional cubes. It
  follows that
  \begin{align*}
    \left( \sum_{k=0}^{2^{nd} - 1} \sum_{r = 1}^{2^d - 1} a_{k, r}
    T_{g_{n, k, r}} \right)(u) & \leq (2^d - 1) 2^{n(d/2-1)} C_P
    \left( \max_{k, r} |a_{k, r}| \right) \sum_{k=0}^{2^{nd} - 1}
    \|Du\|(\operatorname{int} K_{n,k}) \\ & \leq (2^d - 1) 2^{n(d/2-1)} C_P
    \left(\max_{k, r} |a_{k, r}| \right) \|u\|_{BV}.
  \end{align*}
  
  %Let $u \in C^1(I^2)$ a function
  %that we extend to $\R^2$ by zero. For any $k = 0, \dots, 4^n -
  %1$, we have
  %\begin{align*}
  %  a_{3k} T_{g_{n, 3k}}(u) & = a_{3k} \int_{x^{(n, k)}}^{x^{(n, k)} + 2^{-(n+1)}}
  %  \int_{y^{(n, k)}}^{y^{(n, k)} + 2^{-n}} (u(x, y) - u(x + 2^{-(n+1)}, y)) dxdy \\
  %  & = - a_{3k} \int_{x^{(n, k)}}^{x^{(n, k)} + 2^{-(n+1)}}
  %  \int_{y^{(n, k)}}^{y^{(n, k)} + 2^{-n}} \int_0^{2^{-(n+1)}} \frac{\partial u}{\partial x}(x + t, y) dtdxdy \\
  %  & \leq \frac{|a_{3k}|}{2^{n+1}} \int_{K_{n,k}} |\nabla u|
  %  \end{align*}
  %In a very similar fashion, we can show that
  %\[
  %a_{3k + r} T_{n, 3k + r}(u) \leq \frac{|a_{3k+r}|}{2^{n+1}} \int_{K_{n,k}} |\nabla u|
  %\]
  %whenever $k = 0, \dots, 4^n - 1$ and $r \in \{1, 2\}$. It follows
  %that
  %\begin{align}
  %  \left( \sum_{\ell=0}^{3 \cdot 4^n-1} a_\ell T_{g_{n,\ell}} \right)(u) &
  %  \leq \frac{3}{2^{n+1}} \max_{0 \leq \ell \leq 3\cdot 4^n-1} |a_\ell|
  %  \int_{I^2} |\nabla u| \notag \\ & \leq \frac{3}{2^{n+1}}
  %  \left(\max_{0 \leq \ell \leq 3\cdot 4^n-1} |a_\ell| \right)
  %  \|u\|_{BV_\infty} \label{eq:preuveProp5.5}
  %\end{align}
  %Let $B_{BV_\infty(I^2)}$ denote the closed unit ball of
  %$BV_\infty(I^2)$. Theorem~\ref{thm:approx1} and
  %Proposition~\ref{prop:weak*conv} imply that $C^1(I^2) \cap
  %B_{BV_\infty(I^2)}$ is weak* sequentially dense in
  %$B_{BV_\infty(I^2)}$. In view of~\eqref{eq:preuveProp5.5}, this
  %entails that
  %\[
  %\left\| \sum_{\ell=0}^{3 \cdot 4^n-1} a_\ell T_{g_{n,\ell}} \right\| \leq
  %\frac{3}{2^{n+1}} \max_{0 \leq \ell \leq 3\cdot 4^n-1} |a_\ell|
  %\]

  Next, we turn to the lower bound. Let $r \in \{1, \dots, 2^d -
  1\}$. Define
  \[
  u = \sum_{k=0}^{2^{nd}-1} \varepsilon_{k, r} g_{n, k, r}
  \]
  where $\varepsilon_{k, r} \in \{-1,1\}$ are chosen so that
  $\varepsilon_{k, r} a_{k, r} = |a_{k, r}| $. First, we have $\|u\|_1
  = 2^{nd/2}$ and
  \begin{align*}
    \|Du\|(\R^d) & \leq \sum_{k=0}^{2^{nd}-1} \|Dg_{n,k,r}\|(\R^d) \\
    & \leq 2^{nd} 2^{nd/2} \frac{\|Dg_{0,0,r}\|(\R^d)}{2^{n(d-1)}} \\
    & \leq 2^{n(d/2 + 1)} \max \left( \|Dg_{0,0,1}\|(\R^d), \dots,
    \|Dg_{0, 0, 2^{d}-1}\|(\R^d) \right).
  \end{align*}
  Hence, $\|u\|_{BV} \leq 2^{n(d/2+1)} C'$, for some constant $C'$.  As
  the functions $(g_{n, k, r'})$ are pairwise orthogonal in
  $L^2([0, 1]^d)$, we have
  \[
  \left( \sum_{k=0}^{2^{nd}-1} \sum_{r'=1}^{2^d - 1} a_{k, r'}
  T_{g_{n,k, r'}} \right) (u) = \sum_{k=0}^{2^{nd} - 1} |a_{k,r}|.
  \]
  Hence, we infer
  \[
  \sum_{k=0}^{2^{nd} - 1} |a_{k, r}| \leq 2^{n(d/2 + 1)} C' \left\|
  \sum_{k=0}^{2^{nd}-1} \sum_{r' = 1}^{2^d - 1} a_{k, r'} T_{g_{n,k,
      r'}} \right\| .\qedhere
  \]
\end{proof}

\begin{Corollary}
  \label{cor:chargeability}
  Let $f \in C_0([0, 1]^d)$.
  \begin{itemize}
  \item[(A)] The condition
    \[
    \limsup_{n \to \infty} \frac{1}{2^{n(d/2+1)}} \max_{1 \leq r \leq
      2^d - 1} \sum_{k=0}^{2^{nd}-1} |\lambda_{n,k,r}(f)| > 0
    \]
    implies that $f$ is not strongly chargeable.
  \item[(B)] The condition
    \[
    \sum_{n=0}^\infty 2^{n(d/2-1)} \max_{k, r} |\lambda_{n, k, r}(f)|
    < \infty
    \]
    implies that $f$ is strongly chargeable.
  \end{itemize}
\end{Corollary}

\begin{proof}
  Condition (A) implies that the sequence of partial sums
  in~\eqref{eq:decompChargeability} is not Cauchy. Condition (B)
  implies that it is.
\end{proof}

\section{Sample paths of the Brownian sheet}
\label{sec:brownian}

\begin{Empty}
  Let $(\Omega, \calF, \mathbb{P})$ be a probability space and assume it is sufficiently large for the processes below to be
  defined on it.  We recall that the Brownian sheet is a Gaussian
  centered random process $\{W_{t_1, \dots, t_d} : (t_1, \dots, t_d)
  \in [0, 1]^d\}$ with covariance function
  \[
  \Gamma\left( (t_1, \dots, t_d), (t_1', \dots, t_d') \right) =
  \prod_{i=1}^d \min(t_i, t_i').
  \]
  Such a process exists and one may construct it the following way:
  Start from a Gaussian noise $G$ from $L^2([0, 1]^d)$ to a (centered)
  Gaussian space $E$, \textit{i.e.} an isometry from $L^2([0, 1]^d)$ to a
  closed linear subspace $E \subset L^2(\Omega, \calF, \mathbb{P})$
  which contains only centered Gaussian variables. We refer to
  \cite[1.4]{LeGa} for the existence of Gaussian noises. One then sets
  \[
  W_{t_1, \dots, t_d} = G(\ind_{[0, t_1] \times \cdots \times [0,
      t_d]}).
  \]
  This definition leads to the correct covariance function, since
  \[
  \mathbb{E}(W_{t_1, \dots, t_d} W_{t_1', \dots, t_d'}) = \int
  \ind_{[0, t_1] \times \cdots \times [0, t_d]} \ind_{[0, t_1'] \times
    \cdots \times [0, t_d']} = \prod_{i=1}^d \min(t_i, t_i').
  \]
  We let the reader prove that $\Delta_W (K) = G(\ind_K)$ almost surely,
  for a rectangle $K \subset [0, 1]^d$. The Brownian sheet admits a
  continuous modification, by a standard application of Kolmogorov's
  continuity theorem~\cite[Theorem 3.23]{Kall}, so we may suppose that
  $(t_1, \dots, t_d) \mapsto W_{t_1, \dots, t_d}(\omega)$ is an
  element of $C_0([0, 1]^d)$ for all $\omega \in \Omega$.

  As an application of the results from the previous section, we prove
  the following theorem. A generalized result will be obtained in
  section~\ref{sec:fBs} using a more probabilistic argument.
  
  %Note that it is also possible to give a more
  %direct proof of Theorem~\ref{thm:browniansheet} by drawing on the
  %ideas present in the proof of~\cite[Theorem 1.4.5]{Adle}.
\end{Empty}

\begin{Theorem}
  \label{thm:browniansheet}
  The sample paths of the Brownian sheet are almost surely not
  strongly chargeable (for $d \geq 2$).
\end{Theorem}

\begin{proof}
  First we note that $\lambda_{-1}(W) = W_{1, \dots, 1} =
  G(\ind_{[0, 1]^d}) = G(g_{-1})$. Next, using~\eqref{eq:formulaeLambda},
  for $n, k$ we have, almost surely,
  \begin{align*}
    \lambda_{n, k, r}(W) & = 2^{nd/2} \sum_{\ell=0}^{2^d - 1}
    (A_d)_{r, \ell} \Delta_W (K_{n+1, 2^d k + \ell}) \\ & = 2^{nd/2}
    \sum_{\ell=0}^{2^d - 1} (A_d)_{r, \ell} G(\ind_{K_{n+1, 2^d k +
        \ell}}) \\ & = G(g_{n,k,r}),
  \end{align*}
  by~\eqref{eq:defgnkr}. As the sequence of Haar functions $g_{-1},
  g_{0, 0, 1}, \dots$ is orthonormal and $G$ is an isometry, we deduce
  that the random variables $\lambda_{-1}(W), \lambda_{0, 0, 1}(W),
  \dots$ are pairwise uncorrelated and follow the standard Gaussian
  distribution. Since they are jointly Gaussian, we infer that they
  are independent. For each integer $n \geq 0$, define the random
  variables
  \[
  T_n = \frac{1}{2^{nd}} \sum_{k=0}^{2^{nd} - 1} |\lambda_{n , k,
    1}(W)| \text{ and } S_n = 2^{n(d/2-1)} T_n = \frac{1}{2^{n(d/2 +
      1)}} \sum_{k=0}^{2^{nd} - 1} |\lambda_{n , k, 1}(W)|.
  \]
  Each random variable $|\lambda_{n, k, 1}(W)|$ follows a half-normal
  distribution of mean $\sqrt{2 / \pi}$ and variance $1 - 2/\pi$. By
  independence,
  \[
  \left\|T_n - \sqrt{\frac{2}{\pi}} \right\|_2^2 = \operatorname{Var}
  T_n = \frac{1}{2^{nd}} \left(1 - \frac{2}{\pi} \right) \to 0.
  \]
  We infer the existence of a subsequence $(T_{n_k})$ that converges
  almost surely to $\sqrt{2 / \pi}$. Hence $\limsup S_n = \infty$ if
  $d \geq 3$ and $\limsup S_n \geq \sqrt{2/\pi}$ if $d = 2$. In either
  case, $\limsup S_n > 0$ almost surely, and thus we may conclude with
  the help of Corollary~\ref{cor:chargeability}(A).
\end{proof}

\begin{Empty}
We could enhance the previous result by demonstrating that the sample paths of the Brownian sheet are almost surely non-chargeable (though this result is still less comprehensive than Theorem~\ref{thm:counterEx}). This could be accomplished by utilizing the Faber-Schauder basis of the space of charge functionals (as seen in~\ref{e:FaberCH}) rather than that of $SCH([0, 1]^d)$. The approach mirrors the methods employed to establish Theorem~\ref{thm:browniansheet}.
\end{Empty}

\section{Hölder strong charges}

\begin{Empty}
  The 1-dimensional Faber-Schauder basis serves as a wavelet basis for
  the space $C_0([0, 1])$. Consequently, it provides a way to assess
  the regularity of a function in $f \in C_0([0, 1])$ by examining the
  rate at which the coefficients from the Faber-Schauder decomposition
  of $f$ approach $0$.

  The existence of a Faber-Schauder-type basis in the space $SCH([0,
    1]^d)$ supports the idea that some further notions of regularity can be
  formulated for strong charges. In this section, we expound upon a
  theory of Hölder strong charges.

  At first, it is worth noting that if $\mu$ is a strong charge on
  $[0, 1]^d$, and $\alpha = \calS^{-1}(\mu)$ represents its associated
  strong charge functional, then the coefficients $a_{n,k,r}$ of $\alpha$
  in the Faber-Schauder decomposition depend solely on the values on
  $\mu$ over dyadic cubes, since
  \[
  a_{n,k,r} = \alpha(g_{n,k,r}) = 2^{nd/2} \sum_{\ell = 0}^{2^d - 1}
  (A_d)_{r, \ell} \mu(K_{n+1, 2^d k + \ell}).
  \]
  Taking this into consideration, we introduce the following
  definition. Let $\gamma \in \left( \frac{d-1}{d}, 1 \right)$ and
  $\mu$ be a strong charge on $[0, 1]^d$. We say that $\mu$ is
  \emph{$\gamma$-Hölder} whenever there is a constant $C \geq 0$ such
  that $|\mu(K)| \leq C |K|^\gamma$ for all dyadic cubes $K \subset [0, 1]^d$.

  The reader may find it surprising that we have imposed the
  restriction $\gamma > (d-1)/d$ on the Hölder exponent. The reason is
  that there seems to be no meaningful theory for exponents less than or
  equal to $(d-1)/d$. This is foreshadowed in
  Proposition~\ref{analyticalLemma} that provides examples of Hölder
  charges only for exponents greated than $(d-1)/d$.
\end{Empty}

\begin{Empty}[Hölder strong charges in dimension $d=1$]
  In dimension $d=1$, Hölder exponents are allowed to range over
  $\gamma \in (0, 1)$. We claim that a charge $\mu$ is $\gamma$-Hölder
  if and only if its associated continuous function $\varv \colon x
  \mapsto \mu([0, x])$ is $\gamma$-Hölder. This can be proved by elementary methods. In fact, it is well-known that
  a continuous function $\varv \colon [0, 1] \to \R$ is $\gamma$-Hölder
  continuous if and only if there is a constant $C \geq 0$ such that
  \[
  \left| \varv\left( \frac{k+1}{2^n} \right) - \varv\left(
  \frac{k}{2^n} \right) \right| \leq \frac{C}{2^{n\gamma}}
  \]
  for all intergers $n \geq 0$ and $0 \leq k \leq 2^n-1$, see for
  example \cite[Lemma 2.10]{LeGa}.

  In dimension $d \geq 2$, there is no corresponding result. The
  Hölderianity of strong charges is a new regularity notion that has
  no counterpart for functions. In fact, it is possible for a
  $\gamma$-Hölder continuous function in $C_0([0, 1]^d)$ to be
  $\tilde{\gamma}$-Hölder strongly chargeable, where $\tilde{\gamma} >
  \gamma$. The sample paths of the fractional Brownian sheet may
  exhibit this phenomenon, as discussed in section~\ref{sec:fBs}. For
  such functions, adopting the point of view of charges leads to a
  gain in regularity.
\end{Empty}

\begin{Proposition}
  \label{analyticalLemma}
  Let $f \in C_0([0, 1]^d)$ and $\frac{d-1}{d} < \gamma < 1$. Suppose
  that there is a constant $C \geq 0$ such that $|\Delta_f(K)| \leq C
  |K|^\gamma$ for all dyadic cubes $K$. Then $f$ is $\gamma$-Hölder
  strongly chargeable.
\end{Proposition}

\begin{proof}
  By~\eqref{eq:formulaeLambda}, we estimate
  \[
  |\lambda_{n,k,r}(f)| \leq C 2^{nd/2} 2^d \left( \frac{1}{2^{(n+1)d}}\right)^\gamma .
  \]
  Thus,
  \[
  \sum_{n=0}^\infty 2^{n(d/2-1)} \max_{k,r} |\lambda_{n,k,r}(f)| \leq
  C 2^{d(1 - \gamma)} \sum_{n=0}^\infty 2^{n(d-1 -d \gamma)} < \infty,
  \]
  because $d-1 - d \gamma < 0$. By
  Corollary~\ref{cor:chargeability}(B), we conclude that $f$ is
  strongly chargeable. That the strong charge $\Delta_f$ be
  $\gamma$-Hölder now follows from the hypothesis.
\end{proof}

\section{A Kolmogorov-type chargeability theorem for stochastic processes}

\begin{Theorem}[Kolmogorov-type chargeability theorem]
  Let $X$ be a random process indexed on $[0, 1]^d$ with continuous
  sample paths. Let $q > 0, C \geq 0, \delta > 0$ such that
  \begin{equation}
    \label{eq:condKolmo}
  \frac{d-1}{d} < \frac{\delta}{q} \leq 1
  \end{equation}
  and
  \[
  \mathbb{E} \left( |\Delta_X K|^q \right) \leq C |K|^{1 + \delta}
  \]
  for all dyadic cubes $K$. Then a.s., $X$ is $\gamma$-Hölder strongly
  chargeable for any $\frac{d-1}{d} < \gamma < \delta/q$.
\end{Theorem}

\begin{proof}
  It suffices to fix a Hölder exponent $\frac{d-1}{d} < \gamma <
  \frac{\delta}{q}$ and prove the apparently weaker statement that,
  a.s., $X$ is $\gamma$-Hölder strongly chargeable.

  For any integer $p \geq 0$, we let $\calK_p$ be the set of dyadic
  cubes of generation $p$. For such a dyadic cube $K$, we have
  \[
  \mathbb{P} \left( |\Delta_X K| \geq |K|^\gamma \right) \leq
  \frac{1}{|K|^{q\gamma}} \mathbb{E} \left( |\Delta_X K|^q \right)
  \leq C \left( \frac{1}{2^{pd}}\right)^{1 + \delta - q \gamma}.
  \]
  Hence,
  \[
  \mathbb{P} \left(\exists K \in \calK_p :  |\Delta_X K| \geq |K|^\gamma \right)
  \leq C \left( \frac{1}{2^{pd}}\right)^{\delta - q \gamma},
  \]
  as $\# \calK_p = 2^{pd}$. It follows that
  \[
  \sum_{p=0}^\infty \mathbb{P} \left(\exists K \in \calK_p : |\Delta_X
  K| \geq |K|^\gamma \right) < \infty.
  \]
  By the Borel-Cantelli lemma, we have, a.s.,
  \[
  \sup_{p \geq 0} \sup_{K \in \calK_p} \frac{|\Delta_X(K)|}{|K|^\gamma} < \infty.
  \]
  We conclude from Proposition~\ref{analyticalLemma} that the sample paths of $X$
  are almost surely $\gamma$-Hölder chargeable.
\end{proof}

\section{Sample paths of the fractional Brownian sheet}
\label{sec:fBs}

\begin{Empty}
  For any $h \in (0, 1)$, $t, t' \geq 0$, we define
  \[
  \phi^h(t, t') = \frac{|t|^{2h} + |t'|^{2h} - |t - t'|^{2h}}{2}.
  \]
  Let $H_1, \dots, H_d \in (0, 1)$. The \emph{fractional Brownian
    sheet} of Hurst multiparameter $H = (H_1, \dots, H_d)$ is a
  Gaussian centered random process $\{W^H_{t_1, \dots, t_d} : (t_1,
  \dots, t_d) \in [0, 1]^d \}$ of covariance function
  \[
  \Gamma \left( (t_1, \dots, t_d), (t_1', \dots, t_d') \right) =
  \prod_{i=1}^d \phi^{H_i}(t_i, t_i').
  \]
  When $H_1 = \cdots = H_d = 1/2$, we recover the Brownian sheet from the
  next-to-last section. Again, we will suppose that the sample paths of the
  fractional Brownian sheet are continuous. This is possible because
  Kolmogorov's continuity theorem applies to this more general
  case as well. The mean of the Hurst coefficients is written
  \[
  \bar{H} = \frac{H_1 + \cdots + H_d}{d}.
  \]
  This parameter is crucial to determine whether the sample paths of
  the fractional Brownian sheet are chargeable or not. Theorems~\ref{thm:mainFrac} and~\ref{thm:counterEx} describe
  in detail the behavior of these sample paths. We will need to
  compute the variance of the increments in
  Lemma~\ref{lemma:incr}. Though this is well-known, we include this
  calculation for the reader's convenience.
\end{Empty}

\begin{Lemma}[Increments of the fractional Brownian sheet]
  \label{lemma:incr}
  Let $K = \prod_{i=1}^d [a_i, b_i] \subset [0, 1]^d$. Then
  $\Delta_{W^H}(K)$ is a centered Gaussian random variable of variance
  $\prod_{i=1}^d |b_i - a_i|^{2H_i}$. In particular, this variance is
  $|K|^{2\bar{H}}$ when $K$ is a cube.
\end{Lemma}

\begin{proof}
  The random variable $\Delta_{W^H}(K)$ is clearly Gaussian of mean zero, so, we
  need only compute its variance.  We proceed by induction on
  $d$. If $d = 1$ then, clearly,
  \[
  \mathbb{E}\left( (W^{H_1}_{t_1} - W^{H_1}_{t_1'} )^2 \right) =
  \phi^{H_1}(t_1, t_1) - 2 \phi^{H_1}(t_1, t_1') + \phi^{H_1}(t_1',
  t_1') = |t - t'|^{2H_1}.
  \]
  Suppose now that $d \geq 2$ and that the result holds for $d-1$. We
  define a process
  \[
  \widetilde{W}_{t_1, \dots, t_{d-1}} = W^H_{t_1,
    \dots, t_{d-1}, b_d} - W^H_{t_1, \dots, t_{d-1}, a_d}.
  \]
  Its covariance function is
  \begin{align*}
  \lefteqn{\widetilde{\Gamma}\left( (t_1, \dots, t_{d-1}), (t_1',
    \dots, t_{d-1}')\right)} \\ & \qquad = \mathbb{E} W^H_{t_1, \dots,
    t_{d-1}, b_d} W^H_{t_1', \dots, t_{d-1}', b_d} - \mathbb{E}
  W^H_{t_1, \dots, t_{d-1}, b_d} W^H_{t_1', \dots, t_{d-1}', a_d} \\ &
  \qquad \phantom{=} \quad - \mathbb{E} W^H_{t_1, \dots, t_{d-1}, a_d}
  W^H_{t_1', \dots, t_{d-1}', b_d} + \mathbb{E} W^H_{t_1, \dots,
    t_{d-1}, a_d} W^H_{t_1', \dots, t_{d-1}', a_d} \\ & \qquad =
  \left( \prod_{i=1}^{d-1} \phi^{H_i}(t_i, t_i') \right) \left(
  \phi^{H_d}(b_d, b_d) - 2 \phi^{H_d}(a_d, b_d) + \phi^{H_d}(a_d, a_d)
  \right) \\ & \qquad = \left(\prod_{i=1}^{d-1}
  \phi^{H_i}(t_i, t_i') \right)|b_d - a_d|^{2H_d} .
  \end{align*}
  Thus $|b_d - a_d|^{-H_d} \widetilde{W}$ is a fractional Brownian
  sheet of parameter $(H_1, \dots, H_{d-1})$. Observing that
  $\Delta_{W^H}(K) = \Delta_{\widetilde{W}} ([a_1, b_1] \times \cdots
  \times [a_{d-1}, b_{d-1}]) $, one can now conclude.
\end{proof}

\begin{Theorem}
  \label{thm:mainFrac}
  If $\bar{H} > \frac{d-1}{d}$, then a.s., the sample paths of the
  fractional Brownian sheet are $\gamma$-Hölder strongly chargeable
  for any $0 < \gamma < \bar{H}$.
\end{Theorem}

\begin{proof}
  Let $K \subset [0, 1]^d$ be a dyadic cube.  By the Gaussianity of
  the increments $\Delta_{W^H} K$, there is, for each $q > 0$, a
  constant $C_q$ (not depending on $K$) such that
  \[
  \mathbb{E} \left( | \Delta_{W^H} (K) |^q \right) = C_q |K|^{q\bar{H}}.
  \]
  In particular, under the condition that $\bar{H} > (d-1)/d$, the
  Kolmogorov-type chargeability theorem can be applied for $q$ that is
  sufficiently large for
  \[
  q\bar{H} > 1 \text{ and } \bar{H} - \frac{1}{q} > \frac{d-1}{d},
  \]
  and we conclude from this that the sample paths of $W^H$
  are almost surely $\gamma$-Hölder chargeable for any $\gamma <
  \bar{H}$.
\end{proof}

  Next, we will prove that the sample paths of the fractional Brownian
  sheet are almost surely not chargeable whenever $\bar{H} \leq
  \frac{d-1}{d}$. This shows that the
  condition~\eqref{eq:condKolmo} in the chargeability theorem is
  sharp. The proof is based on ideas in
  \cite[Theorem~1.4.5]{Adle}.

  First, we will need the following two lemmas.

\begin{Lemma}
  \label{lemma1}
  Let $X_0, X_1, \dots$ be a sequence of standard normal random
  variables such that for each integer $n$, one has $\lim
  \operatorname{Cov}(X_n, X_{k}) = 0$ as $k \to \infty$. Then
  \[
  \mathbb{P} \left( X_n \geq 1 \text{ infinitely often} \right) = 1.
  \]
\end{Lemma}

\begin{proof}
  Let $(Y_1, Y_2)$ be a Gaussian random vector where
  $\operatorname{Var} Y_1 = \operatorname{Var} Y_2 = 1$ and $\rho =
  \operatorname{Cov}(Y_1, Y_2)$. By the Portmanteau lemma, we have
  \[
  \mathbb{P}(Y_1 \geq 1 \text{ and } Y_2 \geq 1) \xrightarrow[\rho \to
    0]{} \mathbb{P}(Y_1 \geq 1) \mathbb{P}(Y_2 \geq 1).
  \]

  Let $\varepsilon > 0$. By the preceding paragraph, it is possible to
  extract a subsequence $(X_{n_k})$ such that
  \[
  \mathbb{P} \left( X_{n_k} \geq 1 \text{ and } X_{n_\ell} \geq 1 \right) \leq
  (1+ \varepsilon) \mathbb{P}(X_{n_k} \geq 1) \mathbb{P}(X_{n_\ell} \geq 1)
  \]
  for all distinct integers $k, \ell$. It follows that
  \[
  \limsup_{n \to \infty} \frac{\left( \sum_{k = 1}^n \mathbb{P}
    (X_{n_k} \geq 1) \right)^2}{ \sum_{k=1}^n \sum_{\ell=1}^n
    \mathbb{P} (X_{n_k} \geq 1 \text{ and } X_{n_\ell} \geq 1)}
  \geq \frac{1}{1+\varepsilon}.
  \]
  From the Kochen-Stone lemma \cite[Chapter 6, Lemma 4]{FrisGray}, we deduce that
  \[
  \mathbb{P} \left( X_n \geq 1 \text{ i.o.} \right)
  \geq
  \mathbb{P} \left( X_{n_k} \geq 1 \text{ i.o.} \right)
  \geq \frac{1}{1 + \varepsilon}
  \]
  and we finally conclude from the arbitrariness of $\varepsilon$.
\end{proof}

\begin{Lemma}
  \label{lemma2}
  Let $K = \prod_{i=1}^d [a_i, b_i]$ and $K' = \prod_{i=1}^d [a_i',
    b_i']$ be two rectangles in $[0, 1]^d$. Then
  \begin{multline*}
  \operatorname{Cov} \left( \Delta_{W^H} K,
  \Delta_{W^H} K'\right) \\= \frac{1}{2^d}
  \prod_{i=1}^d \left(|b_i' - a_i|^{2H_i} + |b_i - a_i'|^{2H_i} -
    |a_i' - a_i|^{2H_i} - |b_i - b_i'|^{2H_i} \right).
  \end{multline*}
  If $K$ and $K'$ are cubes, then
  \[
  \operatorname{Cov} \left( \frac{\Delta_{W^H} K}{|K|^{\bar{H}}},
  \frac{\Delta_{W^H} K'}{|K'|^{\bar{H}}} \right) = \frac{1}{2^d}
  \prod_{i=1}^d \frac{|b_i' - a_i|^{2H_i} + |b_i - a_i'|^{2H_i} -
    |a_i' - a_i|^{2H_i} - |b_i - b_i'|^{2H_i}}{|b_i - a_i|^{H_i} |b_i'
    - a_i'|^{H_i}}.
  \]
\end{Lemma}

\begin{proof}
  Using the formula for increments, one computes
  \[
    \Delta_{W^H} K \Delta_{W^H} K' = \sum_{(c_i) \in \prod_{i=1}^d
      \{a_i, b_i\}} \sum_{(c_i') \in \prod_{i=1}^d \{a_i', b_i'\}}
    \left( \prod_{i=1}^d (-1)^{\delta_{a_i, c_i} + \delta_{a_i', c_i'}}
      \right) W^H_{c_1, \dots, c_d} W^{H}_{c_1', \dots, c_d'}.
  \]
  Taking expectations on both sides, one finds
  \begin{align}
  \mathbb{E}(\Delta_{W^H} K \Delta_{W^H} K') & = \sum_{(c_i) \in
    \prod_{i=1}^d \{a_i, b_i\}} \sum_{(c_i') \in \prod_{i=1}^d \{a_i',
    b_i'\}} \prod_{i=1}^d (-1)^{\delta_{a_i, c_i} + \delta_{a_i',
      c_i'}} \phi^{H_i}(c_i, c_i') \notag \\ & = \prod_{i=1}^d \left(
  \phi^{H_i}(b_i, b_i') - \phi^{H_i}(a_i, b_i') - \phi^{H_i}(a_i',
  b_i) + \phi^{H_i}(a_i, a_i')\right) \notag \\ & = \prod_{i=1}^d
  \frac{|b_i - a_i|^{2H_i} + |b_i' - a_i|^{2H_i} - |a_i' - a_i|^{2H_i}
    - |b_i' - b_i|^{2H_i}}{2} .\label{eq:covINCR}
  \end{align}
  We recall that the random variables $\Delta_{W^H} K$ and
  $\Delta_{W^H} K'$ are centered. Finally, we get the second equality,
  in case $K$ and $K'$ are cubes, by dividing both sides in the
  equality~\eqref{eq:covINCR} by
  \[
  |K|^{\bar{H}} |K'|^{\bar{H}} = \prod_{i=1}^d |b_i-a_i|^{H_i} |b_i' -
  a_i'|^{H_i} .\qedhere
  \]
\end{proof}

\begin{Theorem}
  \label{thm:counterEx}
  If $\bar{H} \leq \frac{d-1}{d}$, then the sample paths of the
  fractional Brownian sheet are almost surely not chargeable.
\end{Theorem}

\begin{proof}
  To any point $x \in \mathopen{[}0, 1 \mathclose{)}^{d-1}$ and $p \geq
    0$, we associate the dyadic cube $K(x,p)$ in $[0, 1]^d$ defined by
    \[
    K(x, p) = \prod_{i=1}^d \left[ \frac{k_i}{2^p},
      \frac{k_i+1}{2^p}\right], \qquad \text{where } k_i = \lfloor 2^p
    x \rfloor \text{ for } 1 \leq i \leq d - 1 \text{ and } k_d = 0.
    \]
    We note that the sequence of cubes $(K(x,p))$ is decreasing and
    that $(x, 0) \in K(x, p)$.

    We let
    \[
    A = \left\{ (x, \omega) \in \mathopen{[}0, 1 \mathclose{)}^{d-1}
      \times \Omega : \frac{\Delta_{W^H(\omega)} K(x, p)}{|K(x,
          p)|^{\bar{H}}} \geq 1 \text{ for infinitely many } p\right\}
    \]
    where $\Omega$ denotes the underlying sample space on which our
    process is defined. This set is clearly measurable (with respect to the
    product $\sigma$-algebra). For any $x \in \mathopen{[}0,
      1\mathclose{)}^{d-1}$, we define the event $A(x) = \{\omega \in
      \Omega : (x, \omega) \in A)\}$.

    Next, we will apply Lemma~\ref{lemma1} to deduce that $A(x)$ is
    almost certain. To this end, we first notice that the random
    variables $\Delta_{W^H} K(x, p) / |K(x, p)|^{\bar{H}}$ are
    standard normal variables. It then suffices to establish that, for
    any integer $p$,
    \begin{equation} \label{eq:limCov}
    \lim_{q \to \infty} \operatorname{Cov} \left( \frac{\Delta_{W^H}
      K(x,p)}{|K(x,p)|^{\bar{H}}}, \frac{\Delta_{W^H} K(x, q)}{|K(x,
      q)|^{\bar{H}}} \right) = 0.
    \end{equation}
    Suppose $q \geq p$ and write $K(x,p) = \prod_{i=1}^d
    [a_i, b_i]$ and $K(x, q)= \prod_{i=1}^d [a_i', b_i']$. Then the
    covariance in~\eqref{eq:limCov} is given by Lemma~\ref{lemma2}. As
    $K(x, q) \subset K(x,p)$, we can decompose
    \[
    |b_i' - a_i| = |a_i' - a_i| + |b_i' - a_i'| \text{ and }
    |b_i - a_i'| = |b_i - b_i'| + |b_i' - a_i'|.
    \]
    Since the function $x \mapsto x^{2H_i}$ is $\min(2H_i, 1)$-Hölder
    continuous on $[0, 1]$, we infer the existence of a constant $C$
    (depending on $d$, $H_1, \dots, H_d$) such that
    \begin{align*}
    \operatorname{Cov} \left( \frac{\Delta_{W^H}
      K(x,p)}{|K(x,p)|^{\bar{H}}}, \frac{\Delta_{W^H} K(x, q)}{|K(x,
      q)|^{\bar{H}}} \right) & \leq C \prod_{i=1}^d \frac{|b_i' -
      a_i'|^{\min(2H_i, 1)}}{|b_i-a_i|^{H_i} |b_i' - a_i'|^{H_i}} \\ &
    = C \prod_{i=1}^d \frac{|b_i' - a_i'|^{\min(H_i, 1- H_i)}}{|b_i -
      a_i|^{H_i}} \\ & = \frac{C}{|K(x, p)|^{\bar{H}}} \left(\frac{1}{2^q}
    \right)^{\min(H_1, 1 - H_1) + \cdots + \min(H_d, 1 - H_d)}.
    \end{align*}
    This completes the proof of~\eqref{eq:limCov}, from which we can
    assert that $\mathbb{P}(A(x)) = 1$.

    By the measurability of $A$ and the Fubini theorem, it follows
    that almost surely, for almost all $x \in \mathopen{[}0, 1
      \mathclose{)}^{d-1}$, we have $\Delta_{W^H} K(x, p) \geq
      |K(x,p)|^{\bar{H}}$.

    For any integer $n \geq 0$, we are led to consider the (random)
    collection $\calC_n$ of dyadic cubes $K$ of the form $K = K(x,
    p)$, where $p\geq n$ and $\Delta_{W^H} K \geq |K|^{\bar{H}}$. For
    each such cube $K$, we let $\tilde{K}$ denote its ``bottom face''
    $\tilde{K} = \{y \in [0, 1]^{d-1} : (y, 0) \in K\}$, which is
    itself a dyadic cube of $[0, 1]^{d-1}$. Its volume is related to
    that of $K$ by $|\tilde{K}| = |K|^{(d-1)/d}$. By what precedes,
    the cubes $\tilde{K}$, where $K$ ranges over $\calC_n$, cover a
    subset of $[0, 1)^{d-1}$ of full Lebesgue measure. By considering
      only maximal cubes within $\calC_n$, we can extract a finite
      subset $\calD_n \subset \calC_n$ of pairwise disjoint cubes such
      that
      \[
      \sum_{K \in \calD_n} |\tilde{K}| = \sum_{K \in \calD_n} |K|^{\frac{d-1}{d}} \geq \frac{1}{2}.
      \]
      We define $F_n$ to be the dyadic figure $F_n = \bigcup_{K \in
        \calD_n} K$. Then
      \begin{equation}
        \label{eq:contrCBV}
      \Delta_{W^H} F_n = \sum_{K \in \calD_n} \Delta_{W^H} K \geq
      \sum_{K \in \calD_n} |K|^{\bar{H}} \geq \sum_{K \in \calD_n}
      |K|^{\frac{d-1}{d}} \geq \frac{1}{2}.
      \end{equation}
      On top of that, we have $F_n \subset [0, 1]^{d-1} \times [0,
        2^{-n}]$, which ensures that $|F_n| \to 0$. Regarding the
      perimeters, we can estimate
      \[
      \|F_n\| \leq \sum_{K \in \calD_n} \|K\| = 2d \sum_{K \in
        \calD_n} | \tilde{K} | \leq 2d \left|[0, 1]^{d-1}\right| = 2d.
      \]
      Therefore, we have proved that the sequence $(F_n)$
      $BV$-converges to $\emptyset$. The lower
      bound~\eqref{eq:contrCBV} shows that $\Delta_{W^H}$ is not
      continuous with respect to $BV$-convergence and, therefore, cannot
      be extended to a charge. As this happens almost surely, the proof is complete.
\end{proof}

\bibliographystyle{amsplain} \bibliography{phil.bib}

%\printindex

\end{document}